\documentclass{ws-ijbc}
\usepackage{ws-rotating}     
\usepackage{graphicx}
\usepackage{epstopdf}

\usepackage[english]{babel}
\usepackage[utf8]{inputenc}
\usepackage{paralist}
\usepackage[colorinlistoftodos, color=green!40, prependcaption]{todonotes}

\usepackage[pdftex, pdftitle={Article}, pdfauthor={Author}]{hyperref} 
\definecolor{ao}{rgb}{0.0, 0.5, 0.0}
\usepackage{bbm}

\newif\ifFigs       \Figsfalse
\Figstrue

\usepackage{graphicx}

\newcommand{\Includegraphics}[2]%
           {\centering
             \ifFigs\includegraphics[width=#1]{#2}%
             \else  \includegraphics[draft,width=#1]{#2}\fi}

\newcommand{\Includesubgraphics}[3]%
{\begin{minipage}[b]{#1}
    \Includegraphics{\textwidth}{#2}\\
    \centerline{{\footnotesize (#3)}}
  \end{minipage}
}

\usepackage[normalem]{ulem}
\usepackage{xcolor}
\definecolor{blue-violet}{rgb}{0.54,0.17,0.89}
\definecolor{amethyst}{rgb}{0.6,0.4,0.8}
\definecolor{darkviolet}{rgb}{0.58, 0.0, 0.83}
\definecolor{darkgreen}{rgb}{0,.4,0}
\definecolor{mixedgreen}{rgb}{0.3,0.6,00}
\newcommand{\REMOVEva}[1]%
           {{\color{red}\sout{#1}}}

\newcommand{\M}{\mathcal M}

\begin{document}
\title{Doubling bifurcations of invariant closed curves in 3D maps}

\author{Sayanho Biswas
    \footnote{sb21ms164@iiserkol.ac.in}}
    \address{Indian Institute of Science Education $\&$ Research, Kolkata, India}

\author{Soumitro Banerjee
    \footnote{Corresponding author. soumitro@iiserkol.ac.in}}
    \address{Indian Institute of Science Education $\&$ Research, Kolkata, India}

\author{Viktor Avrutin
    \footnote{avrutin.viktor@gmail.com}}
    \address{University of Stuttgart, Germany}

\author{Iryna Sushko
\footnote{sushko@imath.kiev.ua}}
    \address{Institute of Mathematics, National Academy of Sciences of Ukraine, Kyiv, Ukraine}

\maketitle

\begin{abstract}
Numerous studies have reported two types of doubling of invariant closed curves (ICCs) in dynamical systems:
(a) the creation of two disjoint ICCs such that iterations flip between them; and
(b) the creation of a single ICC of double the length of the original ICC.
Various methods have been proposed to analyze such bifurcations. We show that the methods proposed so far are not always adequate for analyzing and predicting these bifurcations. We propose a new method to overcome the difficulties. 
\end{abstract}


\section{Introduction}

An invariant closed curve (ICC) in a Poincar\'e section corresponds to a
two-frequency torus in the phase space of a continuous-time dynamical
system. If the asymptotic dynamics on the ICC is periodic, it is called a
mode-locked ICC or a resonant
torus. In the simplest case, a stable mode-locked ICC contains a pair of cycles, namely a
stable cycle and a saddle, along with the unstable manifold of the
latter converging to the former. Multiple such pairs of cycles may
also be present in a resonant ICC. If the asymptotic dynamics
is quasiperiodic, the curve is called a quasiperiodic ICC or an
ergodic torus, and sometimes also a drift ring. In this case, the
ICC is defined as the closure of any of the quasiperiodic orbits it contains.

An ICC (also called a loop) has been reported to undergo doubling in two ways. In the first type, two disjoint loops develop, and the iterates flip between the two loops. Such bifurcations have been reported in the Mira map \cite{karatetskaia2021shilnikov,gonchenko2021doubling}, mechanical vibro-impacting systems \cite{twoloopmech1,twoloopmech2},
in experiments on molten gallium \cite{TD-moltengallium}, and in
coupled Chua's circuits \cite{2loopchua}. We shall call this kind of
doubling bifurcation a `loop doubling'.

The second type of ICC doubling bifurcation results in a single loop of double length. Such bifurcations have been reported in a map with quadratic nonlinearity \cite{spiegel1983}, in experiments on vibrating strings \cite{strings}, a circuit with hysteresis element
\cite{sekikawa-doubling}, a 3D Lotka-Volterra model
\cite{bischi3Dvolterra}, in a resonant circuit with two saturable
inductors \cite{osinga-tori}, in a bimode laser system
\cite{letellier-laser}, and in the food chain models used in ecological
studies \cite{osipenko}. 
This kind of doubling bifurcation is called `length doubling' (a.k.a. `double covering').

A few approaches have been proposed to analyze the bifurcations of ICCs. In Section~\ref{approaches}, we briefly describe these approaches and show that none of the available procedures is effective in predicting which kind of doubling is going to occur in a system. Moreover, each available approach is applicable to quasiperiodic or mode-locked ICCs, but not both.

In this paper, we propose a new approach based on the topological structure of the tangent space in the doubling direction. Our objective is to develop a general method for predicting the kind of doubling bifurcation that would apply to both resonant and ergodic tori.

\section{State of the art}\label{approaches}

\subsection{\label{sec:2nd:poincare}The second Poincar\'e section}
To analyze bifurcations of quasiperiodic ICCs,
  Banerjee \textit{et al.} \cite{banerjee2012local} proposed placing a `second Poincar\'e section' transverse to the ICC and analyzing the local dynamics around the intersection of the ICC with that plane. Although this approach aptly captures other bifurcations of tori, it cannot distinguish between the two types of ICC doubling because, in the 2nd Poincar\'e section, both are seen as period doubling bifurcation of a fixed point.

\subsection{The sign of the 3rd eigenvalue\label{sec:gardini}}

\begin{figure}[t]
\centering
\includegraphics[width=0.45\textwidth]{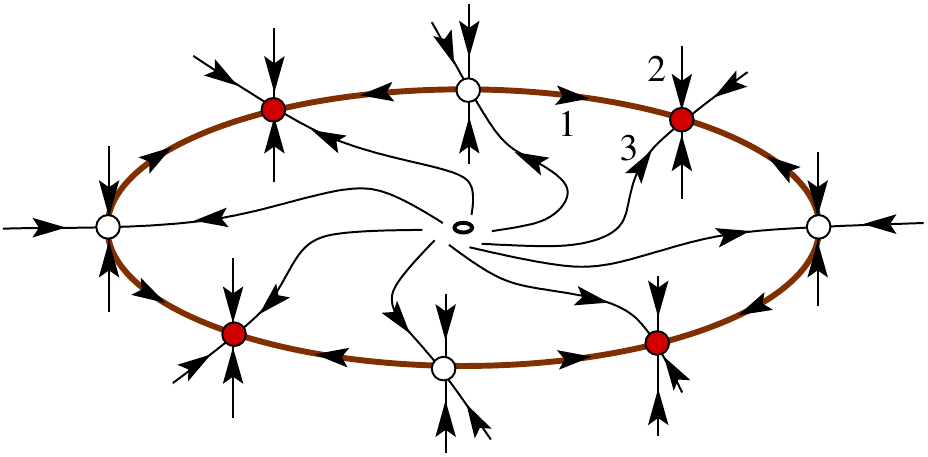}
\caption{Schematic representation of a resonant ICC, showing the eigendirections. The red circles represent the stable cycle, and the white circles represent the saddle cycle. Eigendirection~1 lies along the ICC, eigendirection~2 is the one that bifurcates, i.e., its eigenvalue approaches $-1$. Eigendirection~3 is along the radial direction.\label{schematicloop}}
\end{figure}

This approach proposed by Gardini and Sushko in \cite{gardini2012doubling}, is applicable to mode-locked ICCs. At each point of a stable $p$-cycle belonging to the ICC, there are three eigendirections associated with the eigenvalues of the Jacobian matrix for the $p$th iterate, as schematically shown in Fig.~\ref{schematicloop}. 
In the cited work, the eigenvalue~1 is assumed to be positive (an assumption certainly valid sufficiently close to the Neimark-Sacker bifurcation of a fixed point leading to the appearance of the considered ICC).  In case of a period-doubling of the stable cycle, eigenvalue~2 will necessarily be negative as it goes through $-1$.
In the cited work, it is mentioned that after bifurcation, the manifold that includes ICCs has the form of a normal strip (i.e., a cylinder) in case of loop doubling and a M\"obius strip in case of length doubling. To distinguish between these two cases,
the authors conjectured that the type of doubling depends on the sign of the third eigenvalue: If it is positive, two disjoint ICCs result from the bifurcation, and if it is negative, a single ICC of double length appears.

\subsection{The `Dominant Lyapunov Bundle' method\label{sec:lyapunovbundle}}

This method applies to quasiperiodic ICCs. For a quasiperiodic orbit, the largest Lyapunov exponent is zero, and hence vital information about the orbit is contained in the second-largest Lyapunov exponent, called the `dominant' Lyapunov exponent. The Lyapunov exponents give the rates of exponential divergence or convergence in the direction of the Lyapunov vectors. The collection of Lyapunov vectors corresponding to the dominant Lyapunov exponent computed on the quasiperiodic ICC is called its {\em dominant Lyapunov bundle}.

\begin{figure}[tb]
\centering
\begin{tabular}{ccc}

$\vcenter{\hbox{\includegraphics[width=.23\textwidth]{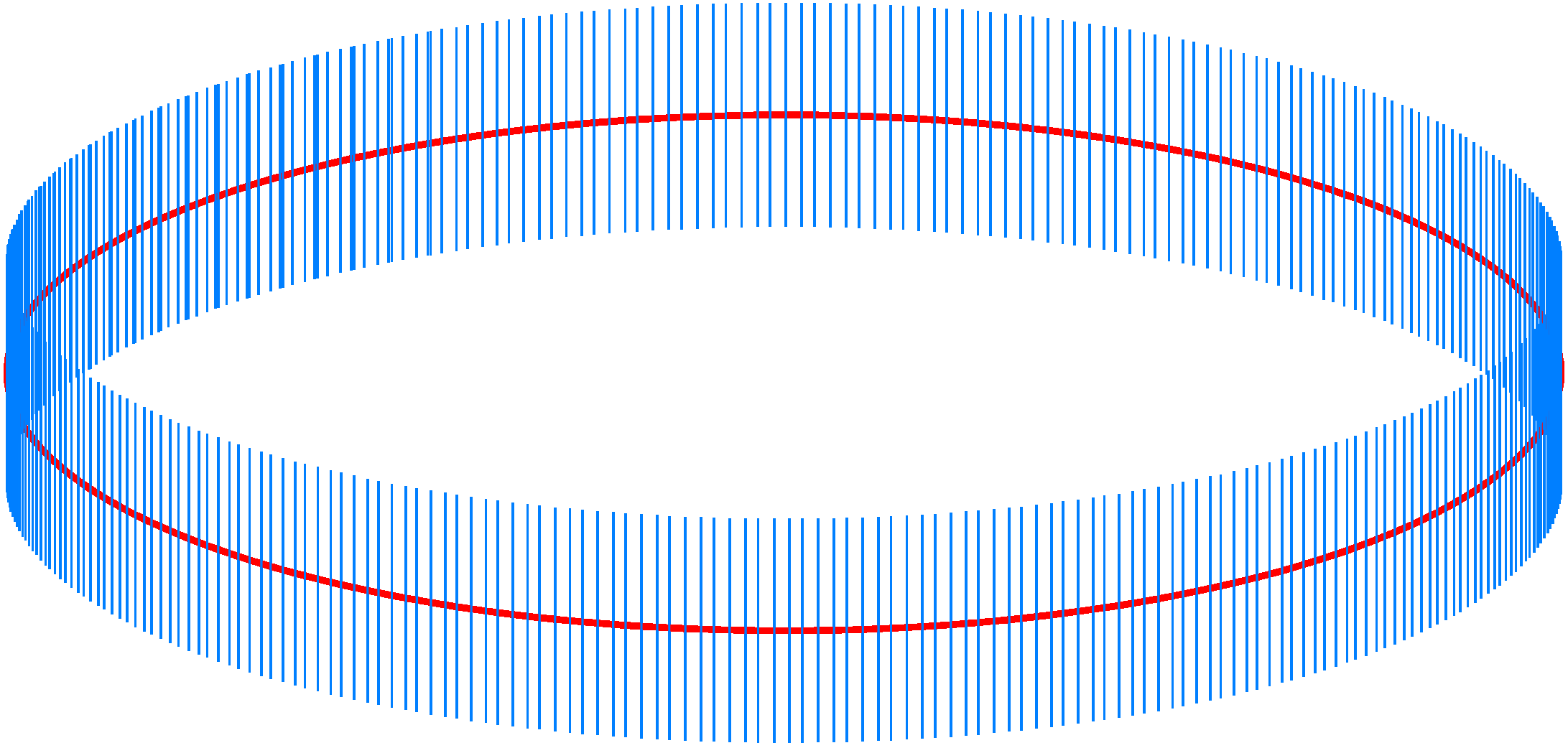}}}$ & ~~$\longrightarrow$~~ & $\vcenter{\hbox{\includegraphics[width=.2\textwidth]{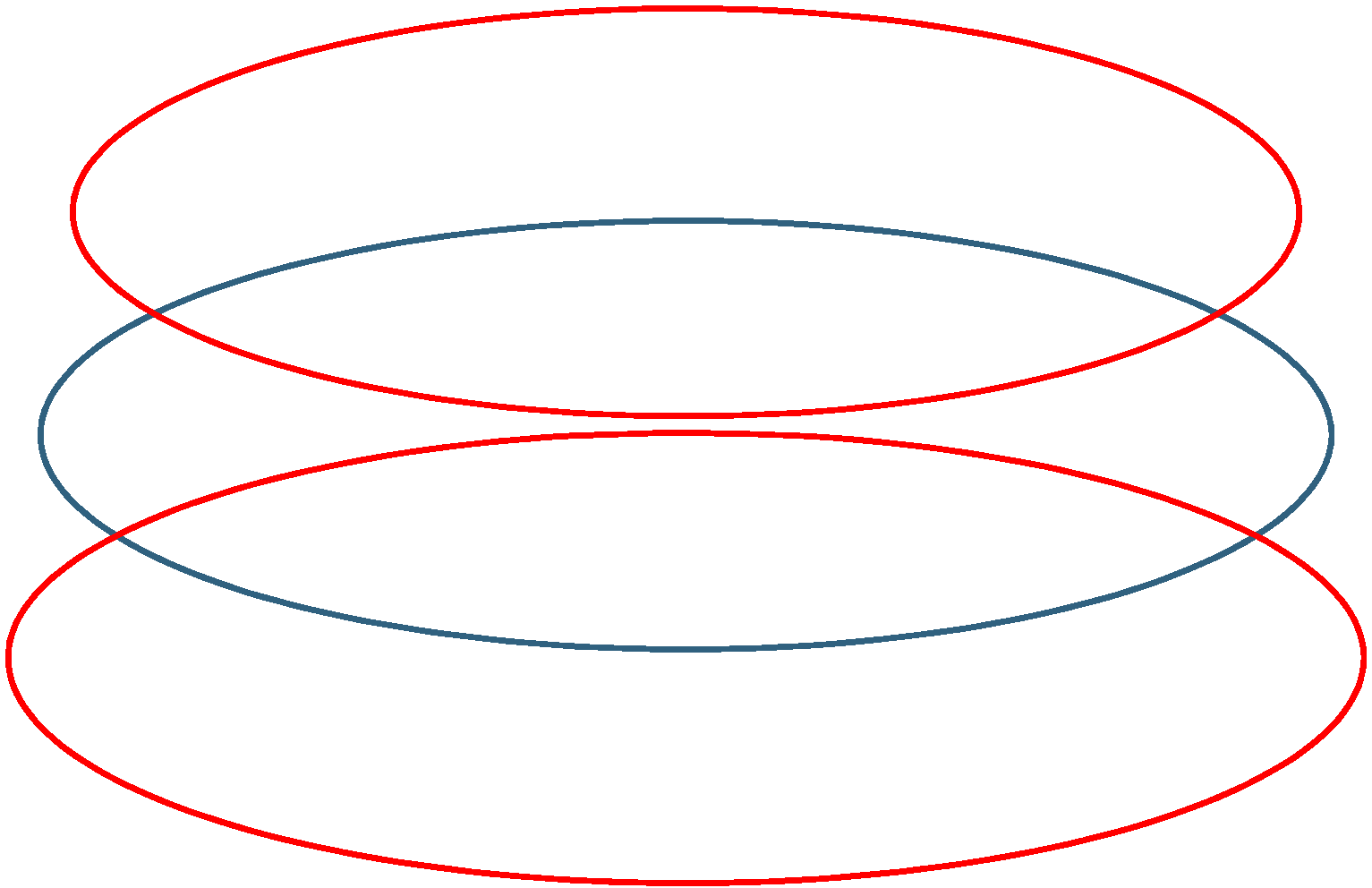}}}$ \\ \multicolumn{3}{c}{\small(a)}\\ \\
$\vcenter{\hbox{\includegraphics[width=.23\textwidth]{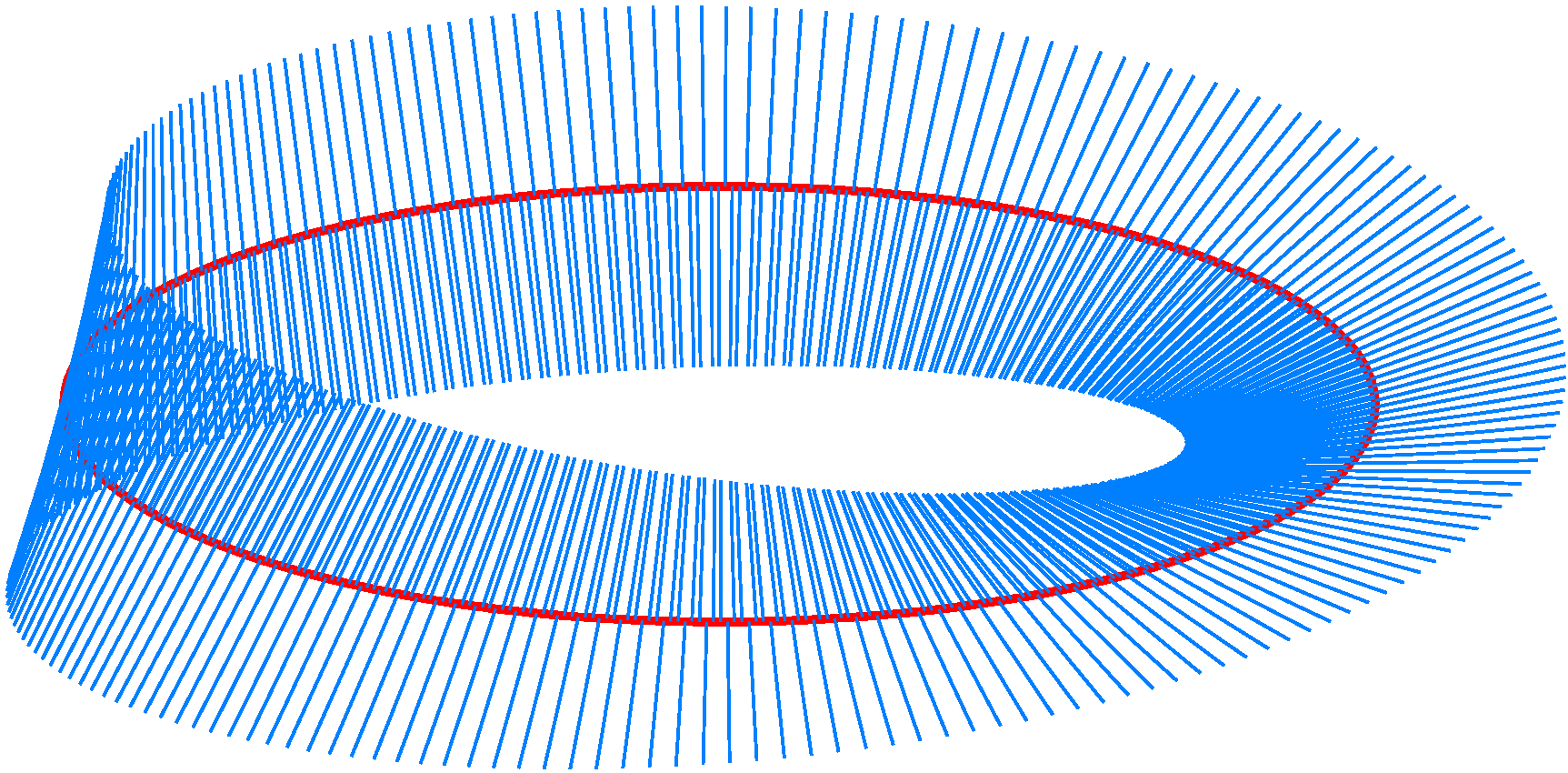}}}$ & ~~$\longrightarrow$~~ &
$\vcenter{\hbox{\includegraphics[width=.2\textwidth]{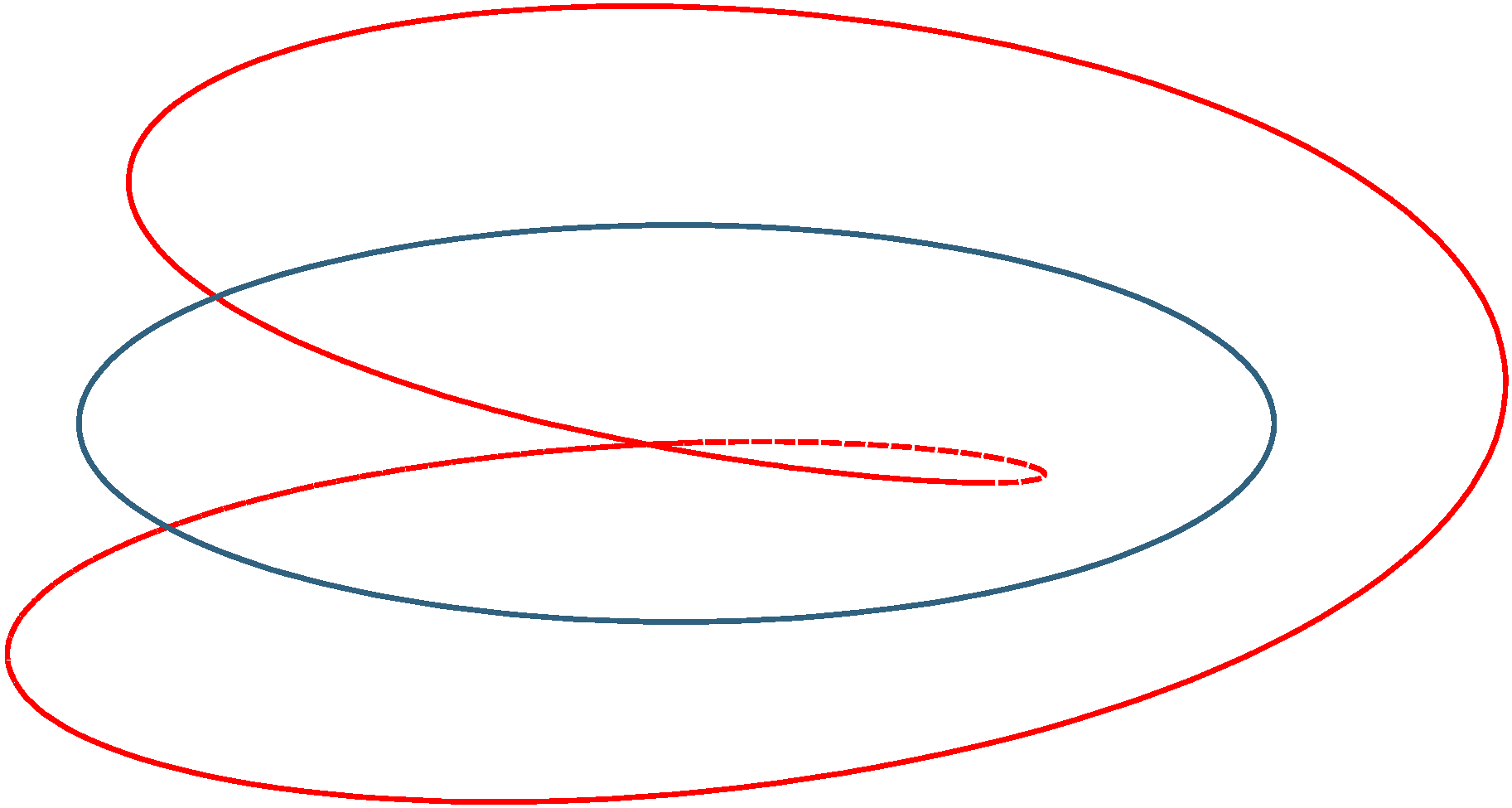}}}$
\\
   \multicolumn{3}{c}{\small(b)}
\end{tabular}

\caption{The topological structures of the dominant Lyapunov bundle close to torus doubling. (a) If it is a cylinder, a loop-doubling bifurcation occurs; (b) If it is a M\"obius strip, a length-doubling bifurcation occurs. The stable orbit is shown in red, and Lyapunov vectors are shown in blue. }
\label{schematic:bifurcations}
\end{figure}

 The method for obtaining Lyapunov bundles is given in \cite{kamiyama2014classification}. The authors showed that, for a loop doubling bifurcation, the topological structure of the dominant Lyapunov bundle is cylindrical (orientation-reversing annulus in their language, see Fig.~\ref{schematic:bifurcations}a). For a length-doubling bifurcation, it takes a M\"obius strip structure (orientation-preserving annulus in their language, see Fig.~\ref{schematic:bifurcations}b). 

 Note that the Lyapunov
bundle method \cite{kamiyama2014classification} is not applicable to resonant ICCs, as the orbits
used to compute the Lyapunov bundles converge to the stable cycle on
the ICC and therefore do not give information on global structures in phase space.

\subsection{The center manifold}

In the work by Gonchenko \textit{et al.} \cite{gonchenko2021doubling}, a more formal foundation for the analysis of bifurcations of
quasiperiodic ICCs has been presented. With the assumption that ICCs must be intersection-free (which is certainly true for invertible maps), they showed that a doubling bifurcation of such a curve is associated with a
2D center manifold and the type of doubling is determined by the topology of this manifold. Specifically, this center manifold has a cylindrical or a M\"obius strip structure, leading to loop and length-doubling bifurcations, respectively.

\subsection{Motivation of the present work}

\begin{compactitem}
\item The existing approaches consider resonant and
quasiperiodic ICCs separately. It is reasonable to assume
that the issue can be treated in a unified way, particularly because quasiperiodic ICCs can often be seen as limiting cases of periodic ICCs (e.g., in period-adding structures).
\item The approaches mentioned in Sections~\ref{sec:2nd:poincare}, \ref{sec:gardini},
and~\ref{sec:lyapunovbundle} have some issues. The secondary Poincar\'e section approach cannot distinguish between loop doubling and length doubling. Later in this paper, we present counter-examples that violate predictions of the `sign of the 3rd eigenvalue' approach. The computation of
Lyapunov bundles is a complicated and computation-intensive procedure, and the method lacks theoretical grounding. 
\item The theoretical foundation provided by 
Gonchenko \textit{et al.} does not provide a way to apply it in the analysis of concrete examples.
\item There has been no investigation on the doubling bifurcation of resonant ICCs involving a saddle-focus connection. This situation needs to be covered under the unified approach.

\end{compactitem}

In this paper, we are interested in developing a practical method
to predict the type of doubling bifurcation
applicable both to resonant and
quasiperiodic ICCs. 

\section{The doubling tangent space}

In this section, we propose a new method for analyzing the doubling of ICCs. 

\subsection{Quasiperiodic ICCs}

     As discussed in \cite{gonchenko2021doubling}, the type of doubling of a quasiperiodic ICC can be predicted based on the topology of the center manifold associated with the eigenvalue $-1$. Unfortunately, this approach is not useful in practical situations, as there is no straightforward method to compute this 2D center manifold. 
     
However, for our purposes, we need neither an analytic representation of this manifold nor its precise approximation, but only the topology. Evidently, near the ICC, this topology is reflected in the topology of a suitable linear tangent space along the ICC, which can be seen as a linear approximation of the relevant 2D center manifold in proximity to the ICC. This tangent space can easily be calculated by evaluating the eigenvectors of the Jacobian matrices for a suitable iterate at the points of the ICC, and it is natural to expect that it preserves the topology of the center manifold.
    
To do this, we first need to calculate the rotation number $\rho$ of the ICC. If the ICC is quasiperiodic, $\rho$ is an irrational number. We choose a sufficiently good rational approximation (one could use a Diophantine approximation) $\rho \approx q/p$.  We evaluate the eigenvalues of the Jacobian matrix $J_{f^p}$ of the $p$th composition of the map $f$ at the points of the ICC. Since we consider the ICC sufficiently close to a doubling bifurcation, each of these matrices has an eigenvalue close to $-1$. The union of the eigenvectors corresponding to these eigenvalues computed at points on the ICC approximates the tangent eigenspace. 
This tangent eigenspace reflects the topology of the 2D center manifold, and hence can be used instead of the center manifold to predict the type of doubling bifurcation.

Next, note that the concept of a center manifold is applicable only at the bifurcation point. 
However, for a smooth system sufficiently close to
the bifurcation, there exists a slow manifold of the
same topology. The structure of this
manifold changes to a small extent for small changes of a parameter, and, like the center manifold at the bifurcation point, it can be approximated by
the corresponding linear eigenspace.
Therefore, the topology of the tangent eigenspace before the bifurcation predicts the topology of the tangent eigenspace at the bifurcation which, in turn, reflects the topology of the center manifold at the bifurcation point. 
Hence, the topology of the tangent eigenspace before the bifurcation predicts the kind of doubling bifurcation the ICC would eventually undergo.

\subsection{Resonant ICCs}

\medskip
\subsubsection*{Setup}

Let us first fix our setup.  We consider a resonant ICC in a three-dimensional map $f$, consisting of a stable period-$p$ node, a period-$p$ saddle, and the unstable manifold of the saddle converging
to the node. To describe the dynamics on the ICC, we consider the eigenvalues and eigenvectors of the Jacobian $J_{f^p}$ of the $p$-th iterate of the map, for which the cycle points are fixed points.

Each cycle has three eigendirections determined by $J_{f^p}$, as
illustrated in Fig.~\ref{schematicloop}.  For simplicity, we denote these directions as
radial, tangential, and transversal.  This notation is borrowed from
the configuration near a Neimark-Sacker bifurcation, where the ICC can
be approximated by a circle in a plane.  In this context, the radial
direction corresponds to the radius of the circle, the tangential
direction aligns with the tangent to the circle in the plane, and the
transverse direction is orthogonal to this plane.

Note that in publications related to resonant ICCs,
several mistakes have been made by mixing the
eigenvalues of the Jacobians of $f$ and $f^p$.
Throughout this paper, we use the Jacobian of $f^p$.

\subsubsection*{Assumptions}

Regarding the eigenvalues of $J_{f^p}$ associated with the three
directions described above, we make the following assumptions.
\begin{compactenum}[{A}1.]
\item 
  The eigenvalue associated with the radial direction is positive and less than one for both cycles (consistent with dynamics soon after a
Neimark-Sacker bifurcation).
  
\item 
  The eigenvalue associated with the tangential direction is positive for both cycles: it is less than one for the stable node and greater than one for the saddle (as expected after a fold bifurcation determining the boundaries of the related periodicity tongue issuing from a Neimark-Sacker bifurcation boundary).

\item 
  The eigenvalue associated with the transverse direction is negative for both cycles. Initially, it lies between $-1$ and $0$, but under parameter variation, these eigenvalues approach and eventually cross $-1$. Consequently, the period-doubling of the ICC occurs in the transversal direction. 
\end{compactenum}

The saddle and node belonging to the resonant ICC generally do not undergo doubling bifurcations at the same parameter value. Therefore, the doubling of a resonant ICC occurs not through a single bifurcation but via a sequence of two bifurcations. 
Thus, two distinct doubling scenarios arise, depending on whether the saddle or the node doubles first. 

The two scenarios are independent of the type of doubling exhibited by the ICC (loop or length doubling). Moreover, the doubling bifurcation of one cycle does not imply that the other cycle must also undergo a doubling bifurcation.
Therefore, we introduce an additional assumption:
\begin{compactenum}[{A}4.]
\item
  The doubling bifurcations of the saddle and node belonging to the resonant ICC occur at parameter values sufficiently close to each other, and no other bifurcations involving these cycles arise between these values.
\end{compactenum}

\subsubsection*{Main result}
Under assumptions A1-A4, we develop a way of predicting the type of doubling exhibited by a resonant ICC. Note that at the point of bifurcation of the saddle or node, the center manifolds are one-dimensional curves that are transversal to the ICC. These 1D center manifolds do not yield conclusive information regarding the type of doubling. 

Recall that, in the case of a quasiperiodic ICC, the type of doubling it undergoes can be inferred from the topology of a two-dimensional center manifold. We proceed to define a similar two-dimensional attracting manifold $\M$ with the
following properties:
\begin{compactenum}[1.]
\item     
  The ICC is contained within $\M$.
\item At the parameter values corresponding to the doubling bifurcations,
  the one-dimensional center manifolds are also contained within $\M$.
\item 
 At other points along the ICC, $\M$ is tangent to the eigenvectors that, at the points of the stable cycle and at the bifurcation parameter value, correspond to the eigenvalue $-1$.
\end{compactenum}

\begin{figure}[t]
    \centering
    \includegraphics[width=0.3\textwidth]{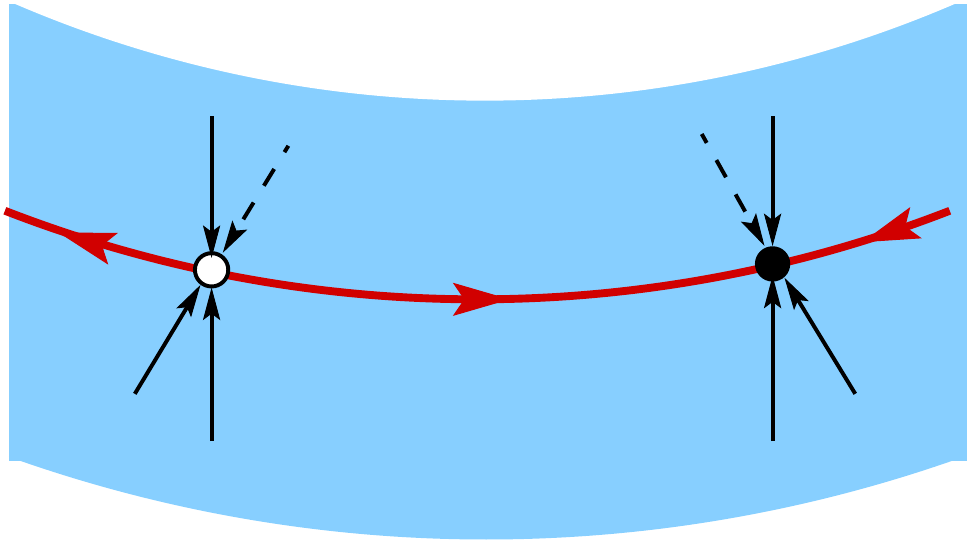}
    \parbox[b][9em][c]{3em}{~~\Large $\longrightarrow$~~} 
    \includegraphics[width=0.25\textwidth]{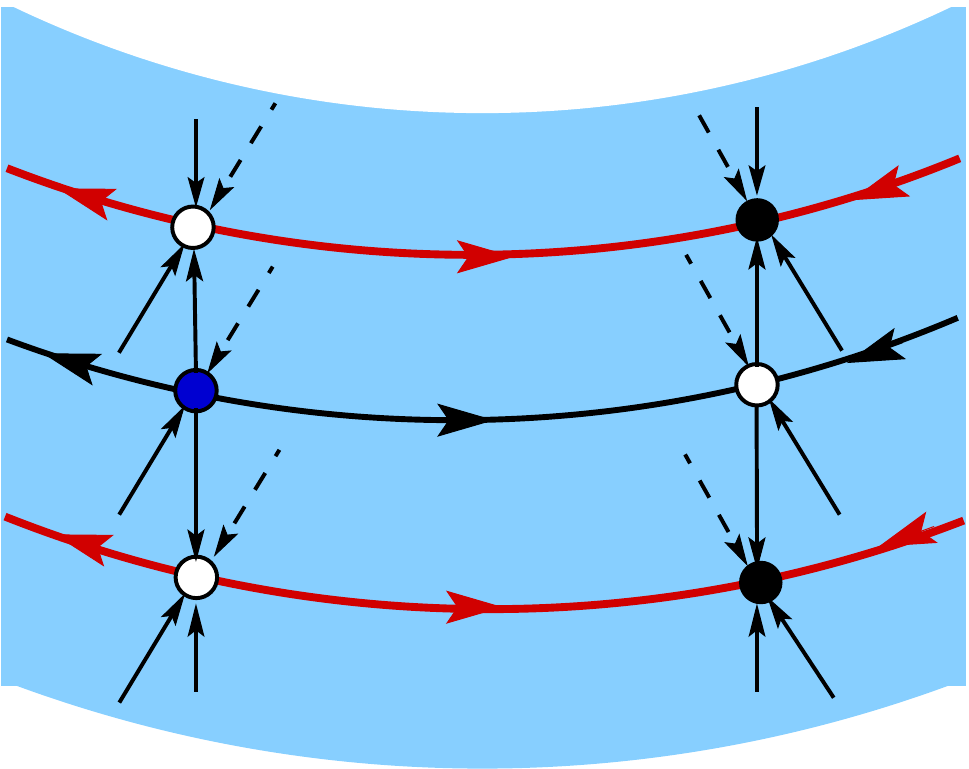}
    \caption{Schematic picture of the attracting manifold $\M$, showing the local structure before and after the doubling of the ICC.} \label{sn}
\end{figure}

We claim that such an attracting surface exists because of the strong contraction in the radial direction (see Fig.~\ref{sn}). If assumption A4 is satisfied, $\M$ would persist throughout the parameter interval between the doubling bifurcations of both cycles.  The key statement of our theory is as follows.\\

  \emph{If the two-dimensional manifold $\M$ exists,  its topology determines the type of doubling: either cylindrical or M\"obius-strip-like --- corresponding respectively to
  loop- or length-doubling of the resonant ICC.}

\subsubsection*{Computation}

To determine the type of doubling a resonant ICC
undergoes, we must determine the topology of $\M$.  As with quasiperiodic ICCs, this topology is reflected --- ultimately via the
Hartman-Grobman theorem --- in the topology of the corresponding tangent space.  When analyzing resonant ICCs, we apply the
same computational procedure used for quasiperiodic ICCs. We compute the eigenvectors of the matrices $J_{f^p}$
associated with the doubling (transversal) direction at points along
the ICC --- both the cycle points and the points on the saddle's
unstable manifold.  In contrast to the quasiperiodic case, we do not
need to approximate the rotation number $\rho$, as the period $p$ is
known a priori.

\subsubsection*{Dynamics inside $\M$}

We now return to the two doubling scenarios, which depend on whether
the saddle or the node doubles first. Due to the
contracting nature of the radial direction, these scenarios unfold within $\M$.

\medskip

\noindent {\bf If the saddle doubles first:} 

\begin{figure*}[t]
\centering
\fbox{\includegraphics[height=1.1in, width=1.5in]{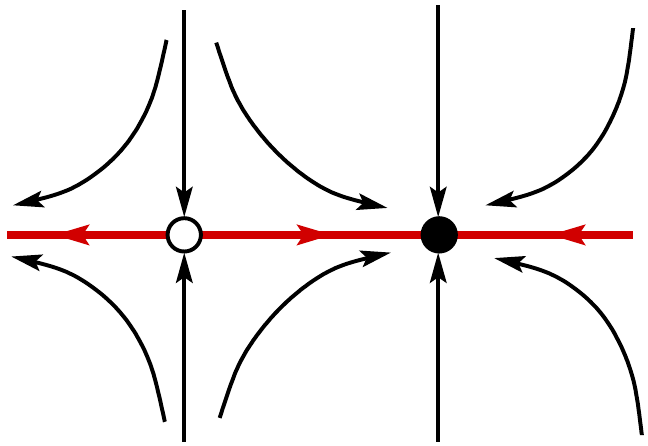}{\footnotesize(a)}}
 \parbox[b][7em][c]{3em}{~~\Large$\longrightarrow$~~} 
     \fbox{\includegraphics[height=1.1in, width=1.5in]{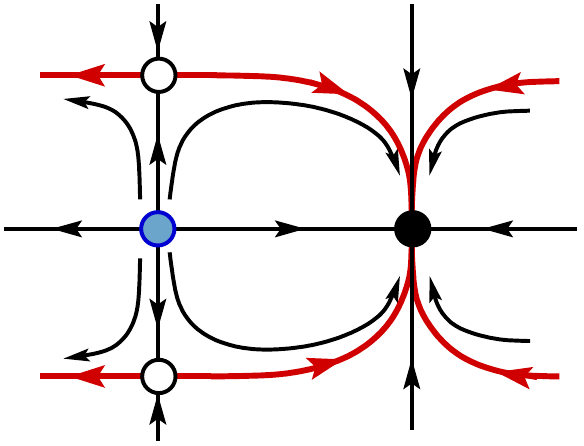}{\footnotesize(b)}}   
     \parbox[b][7em][c]{3em}{~~\Large $\longrightarrow$~~}   \fbox{\includegraphics[height=1.1in, width=1.5in]
      {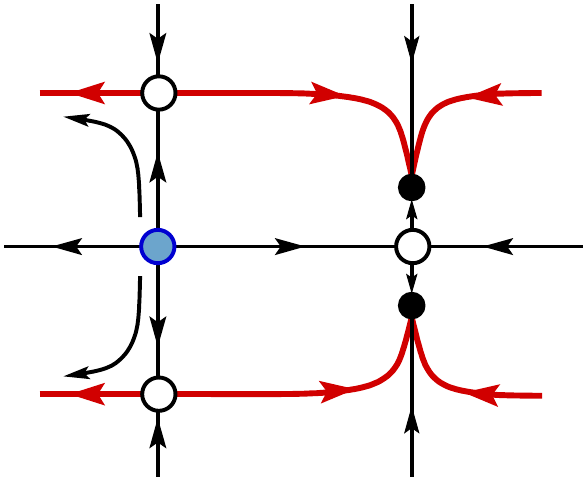}{\footnotesize(c)}}  \\
      \medskip
 {\Large $\downarrow$} \hspace{4.5in} {\Large $\downarrow$}\\ \medskip
\fbox{\includegraphics[height=1.1in, width=1.5in]{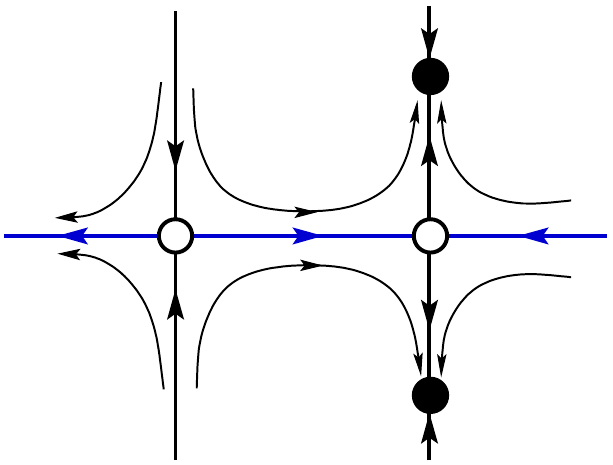}{\footnotesize(d)}}
  \parbox[b][7em][c]{3em}{~~\Large $\longrightarrow$~~} 
\fbox{\includegraphics[height=1.1in, width=1.5in]{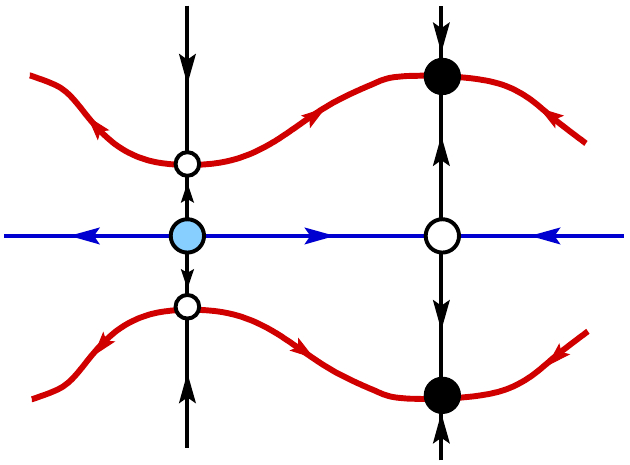}{\footnotesize(e)}}   
     \parbox[b][7em][c]{3em}{~~\Large $\longrightarrow$~~} 
    \fbox{\includegraphics[height=1.1in, width=1.5in]{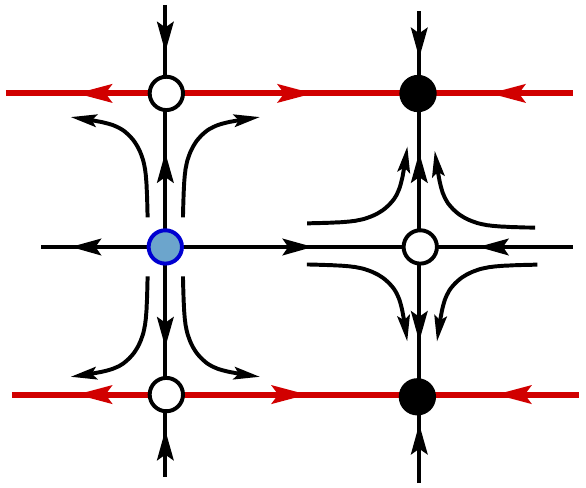}{\footnotesize(f)}}

      \medskip
    \caption{\label{resonant-bif} Schematic picture of dynamics on $\M$. (a) The saddle-node connection; (b) The saddle has doubled but the node has not; (c) The node has just doubled; (d) The node has doubled but the saddle has not; (e) The saddle has also doubled; (f) The final situation. The system takes either (a) $\rightarrow$ (b) $\rightarrow$ (c) $\rightarrow$ (f) route, or (a) $\rightarrow$ (d) $\rightarrow$ (e)  $\rightarrow$ (f) route.}
\end{figure*}

The process starts with an ICC with a saddle-node connection; see Fig.~\ref{resonant-bif}(a).

In the first step, the original saddle becomes a repelling node. In its vicinity, a period-$2p$ saddle emerges. As a result, an attracting bilayered ICC arises, consisting of a period-$p$ node, a period-$2p$ saddle, and its unstable manifold (see Fig.~\ref{resonant-bif}(b)). At this stage, the original period-$p$ ICC no longer exists in the usual sense, as the repelling node has a two-dimensional unstable manifold (a subset of $\M$), while the stable node possesses a stable two-dimensional manifold (a subset of $\M$ as well). 
  
In the next stage, the original node becomes a saddle. In its vicinity, a period-$2p$ node appears. At the bifurcation, the branches of the unstable manifold of the period-$2p$ saddle are tangent to the eigenvectors associated with the eigenvalue $-1$ (Fig.~\ref{resonant-bif}(b)).

After the bifurcation, a period-$2p$ ICC is present, consisting of a stable period-$2p$ node, a period-$2p$ saddle, and its unstable manifold (Fig.~\ref{resonant-bif}(e)). Note that, immediately after bifurcation, the manifold branches approach the saddle from the same side (see Fig.~~\ref{resonant-bif}(c)). In addition, the unstable period-$p$ ICC is restored, now composed of the unstable period-$p$
  node, a period-$p$ saddle and its stable manifold.

\medskip

\noindent {\bf If the node doubles first:} 

The sequence of events starts again with a stable period-$p$ ICC as shown in Fig.~\ref{resonant-bif}(a).

In the first stage, the original node becomes a saddle. In its vicinity, a stable period-$2p$ node emerges (see Fig.~\ref{resonant-bif}(d)). As a result, a repelling period-$p$ ICC is formed, containing two period-$p$ saddles---one with an attracting transversal direction and the other with a repelling one. The unstable manifold of the first coincides with the stable manifold of the second, forming the ICC. At this stage, the period-$2p$ attracting node is not yet involved in any ICC.

In the second stage, the original saddle --- previously with an attracting transverse direction --- becomes a repelling node. In its vicinity, a period-$2p$ saddle emerges. After bifurcation, an attracting period-$2p$ ICC is established (Fig.~\ref{resonant-bif}(f)).  The unstable period-$p$ ICC exists as well (see Fig.~\ref{resonant-bif}(e)), which includes the period-$p$ repelling node, the period-$p$ saddle, and an unstable manifold of the former connecting to the latter.

\subsection{A note on cycles of even period}
\begin{figure}
    \centering
    \includegraphics[width=0.23\textwidth]{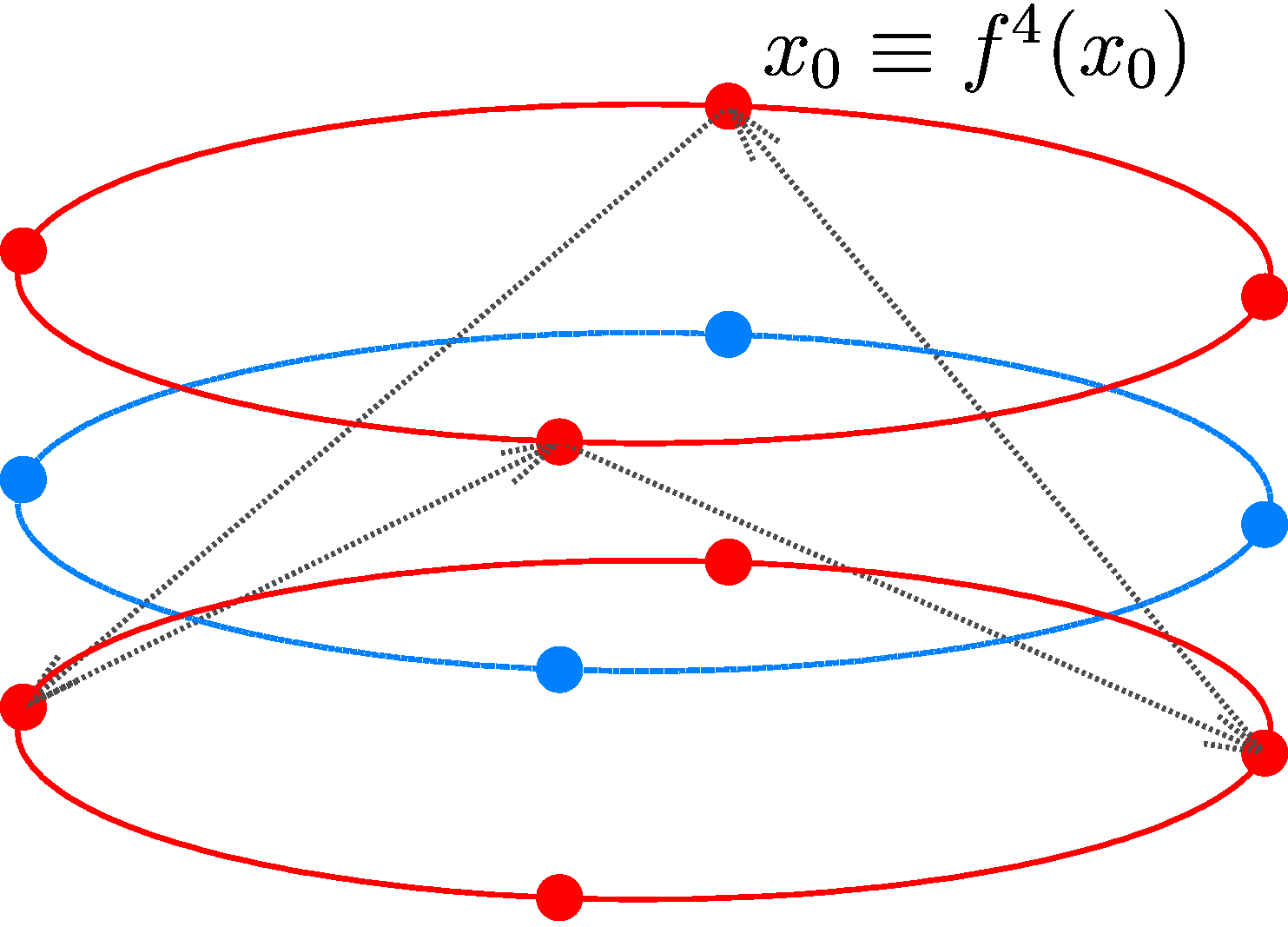}\hspace{0.3in}
     \includegraphics[width=0.24\textwidth]{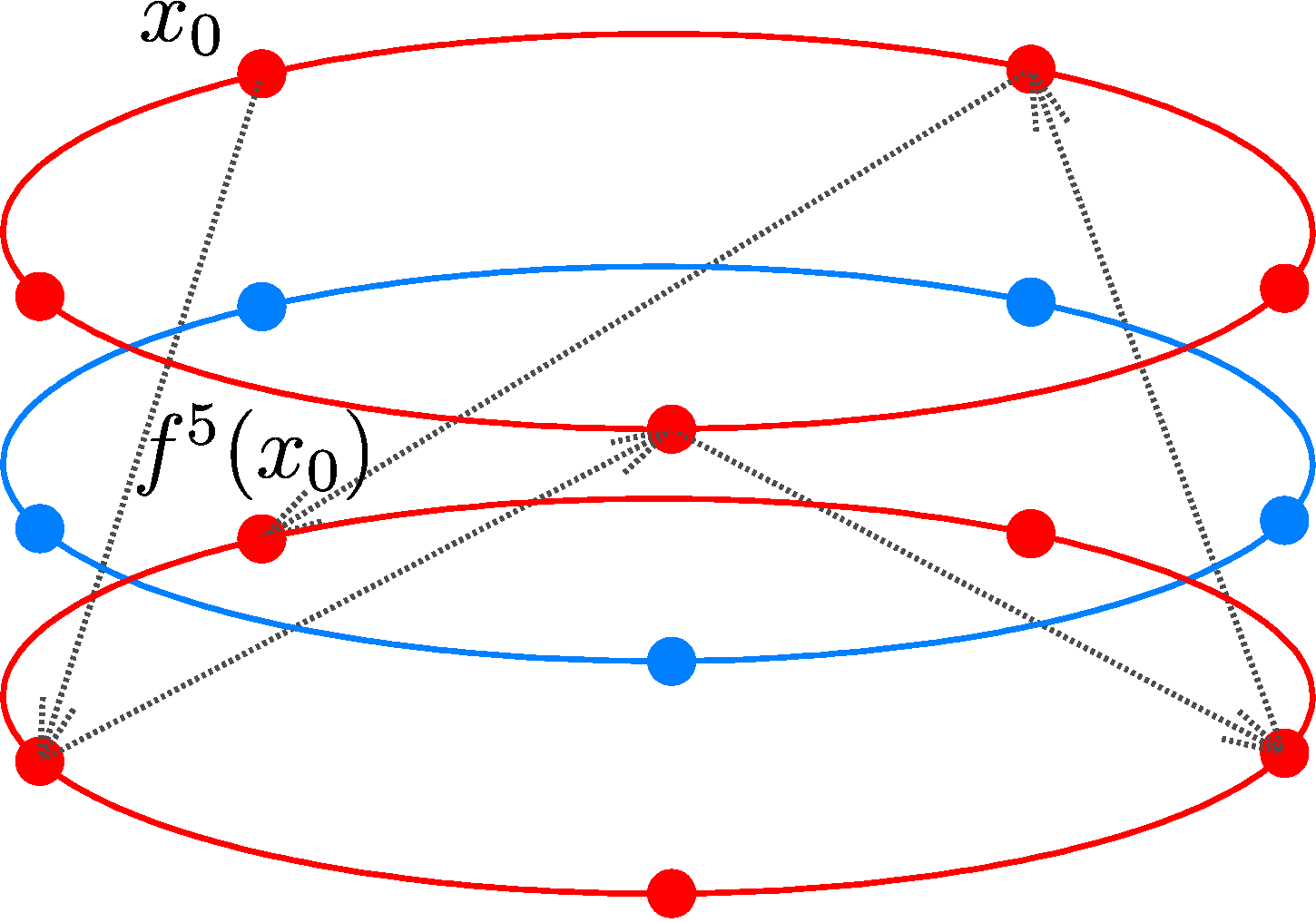}
    \caption{    
    Schematic representation of two ICCs and a $2p$-cycle (in red) after a loop doubling of a resonant ICC (in blue). For each curve, only one of the two cycles is shown to avoid overloading the figure. Here, $p$ is (a) even (not possible); (b) odd (possible). The iterates of a cycle point 
    $x_0$ are indicated by gray arrows.
   }\label{fig:odd:even}
\end{figure}

As explained above, the type of doubling can be predicted based on the topology of the tangent doubling eigenspace. However, it is straightforward to show that for a cycle of even period $p$, the only possible topology of this eigenspace is that of a M\"obius strip. To demonstrate this, suppose that the doubling eigenspace has a cylindrical topology. In this case, the second iterate $f^2$ of the map after the doubling bifurcation must have two separate ICCs, each with an attracting $p$-cycle, together forming an attracting $2p$-cycle for $f$. This means that the iterates of a point $x_0$ of this $2p$-cycle under $f$ jump from one ICC to the other at each step. However, this is not possible, because iterates that jump in this manner return to the point $x_0$ after $p$ steps (i.e., $x_0=f^p(x_0)$), rather than $2p$ (see Fig.~\ref{fig:odd:even}(a) where $p=4$). This implies that the topology of the tangent eigenspace cannot be cylindrical and must instead be a M\"obius strip. For odd periods, this problem does not occur and the cylindrical topology is possible (see Fig.~\ref{fig:odd:even}(b) with $p=5$).

\subsection{Practical issues}

At this stage, some practical issues are worth mentioning.

\medskip
\begin{compactitem}
\setlength{\itemsep}{0.6em}
\item Note that computing the eigenvalues
 of $J_{f^p}$ for large $p$ may be an ill-conditioned task.
 
\item The topology (cylinder or M\"obius) of the tangent eigenspace may not always be clear from its appearance, since the eigenspace can be arbitrarily twisted.  Still, determining this topology is a simple task, as one can add the angles between the eigenvectors at the neighboring points of the ICC and the sum ($\bmod\: 2\pi$) converges either to zero or to $\pi$. If the sum yields zero, it has a cylindrical topology and if it yields $\pi$, it is M\"obius.   

\item For quasiperiodic ICCs, both Lyapunov bundles and the tangent eigenspaces are approximations of the 2D center manifold. However, the tangent eigenspace is computationally much easier to obtain than the dominant Lyapunov bundles. 

\item 
For resonant ICCs, the computation procedure is more complicated than
for quasiperiodic ones, as it includes computation of the ICCs itself, which, 
in this case, cannot be obtained by forward iterations starting from the points of the $p$-cycle. One has to detect the saddle cycle and then compute the
unstable manifold of the saddle. When dealing with resonant ICCs of a sufficiently high period, this step might be unnecessary if the topology of the tangent eigenspace
is evident from the union of doubling eigenvectors at the points of the stable
cycle.

\end{compactitem}

\section{The considered systems }

In this section, we introduce the systems used to illustrate the above theoretical results.

\subsection{The Mira map}
We consider the Mira map \cite{MirGarBarCat96}, defined by
\begin{equation}
\label{eq:mira:map}
   \begin{split}
    x_{n+1} &= y_n\\
    y_{n+1} &= z_n\\
    z_{n+1} &= b x_n + c y_n + a z_n - x_n^2
   \end{split} 
\end{equation}

The map has a trivial fixed point
at $x=y=z=0$. The trace, determinant, and second trace of the
Jacobian, evaluated at this fixed point, are $\tau=a$, $\delta=b$, and
$\sigma=-c$, respectively.  According to the generic Neimark-Sacker bifurcation
condition for fixed points in 3D maps, namely $\delta^2 -
\tau\delta+\sigma=1$, we conclude that this fixed point undergoes a
Neimark-Sacker bifurcation at the set
\begin{equation}
  \phi_{\text{\tiny NS}} = \{(a,b,c) \mid b^2 - ab - c = 1\}
\end{equation}

\begin{figure*}[t!]\centering
  \Includesubgraphics{.4\textwidth}{fig_Mira_2Da}{a}\hspace{0.3in}
  \Includesubgraphics{.4\textwidth}{fig_Mira_2Db}{b}\\
  \medskip
  \Includesubgraphics{.4\textwidth}{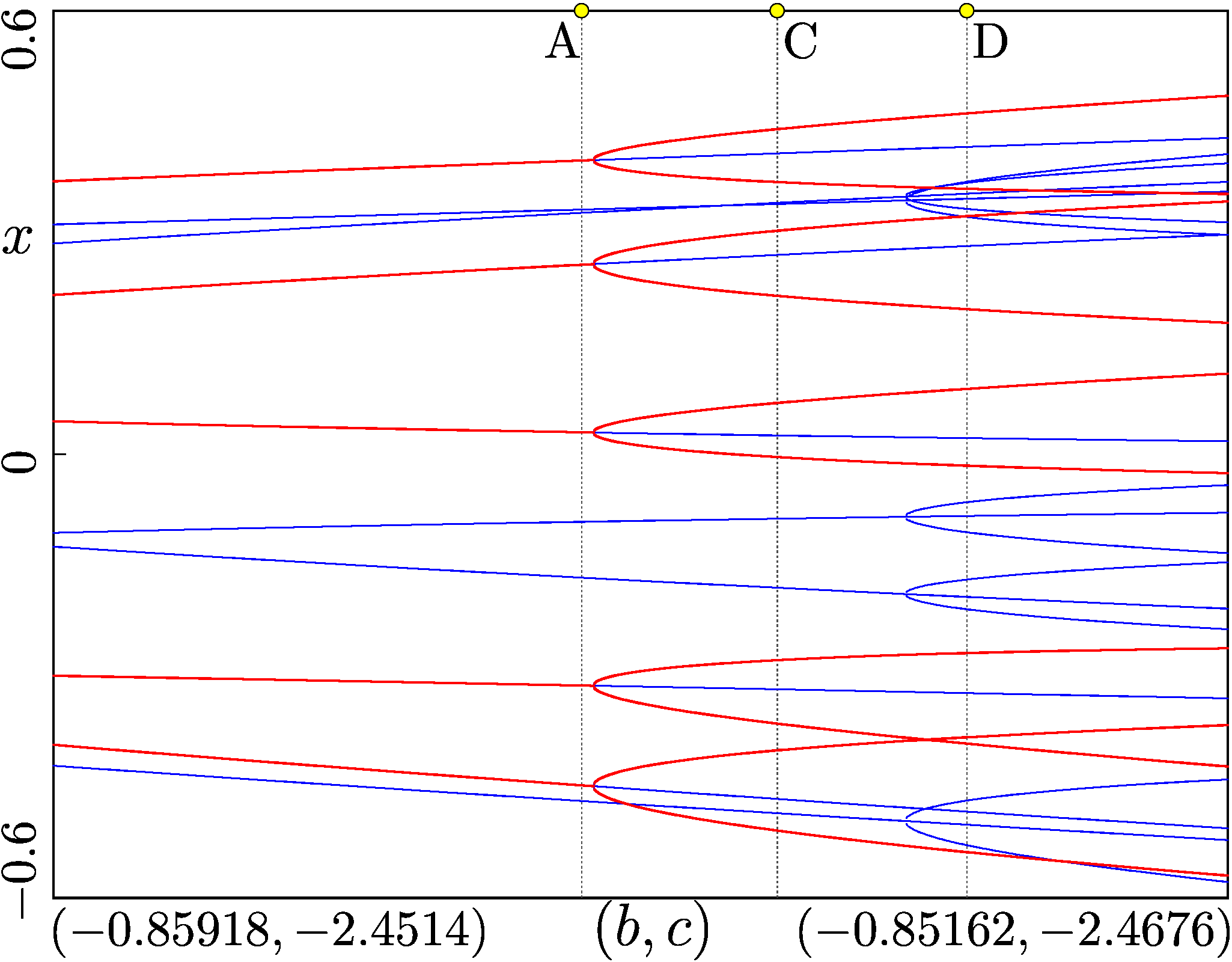}{c}\hspace{0.3in}
  \Includesubgraphics{.4\textwidth}{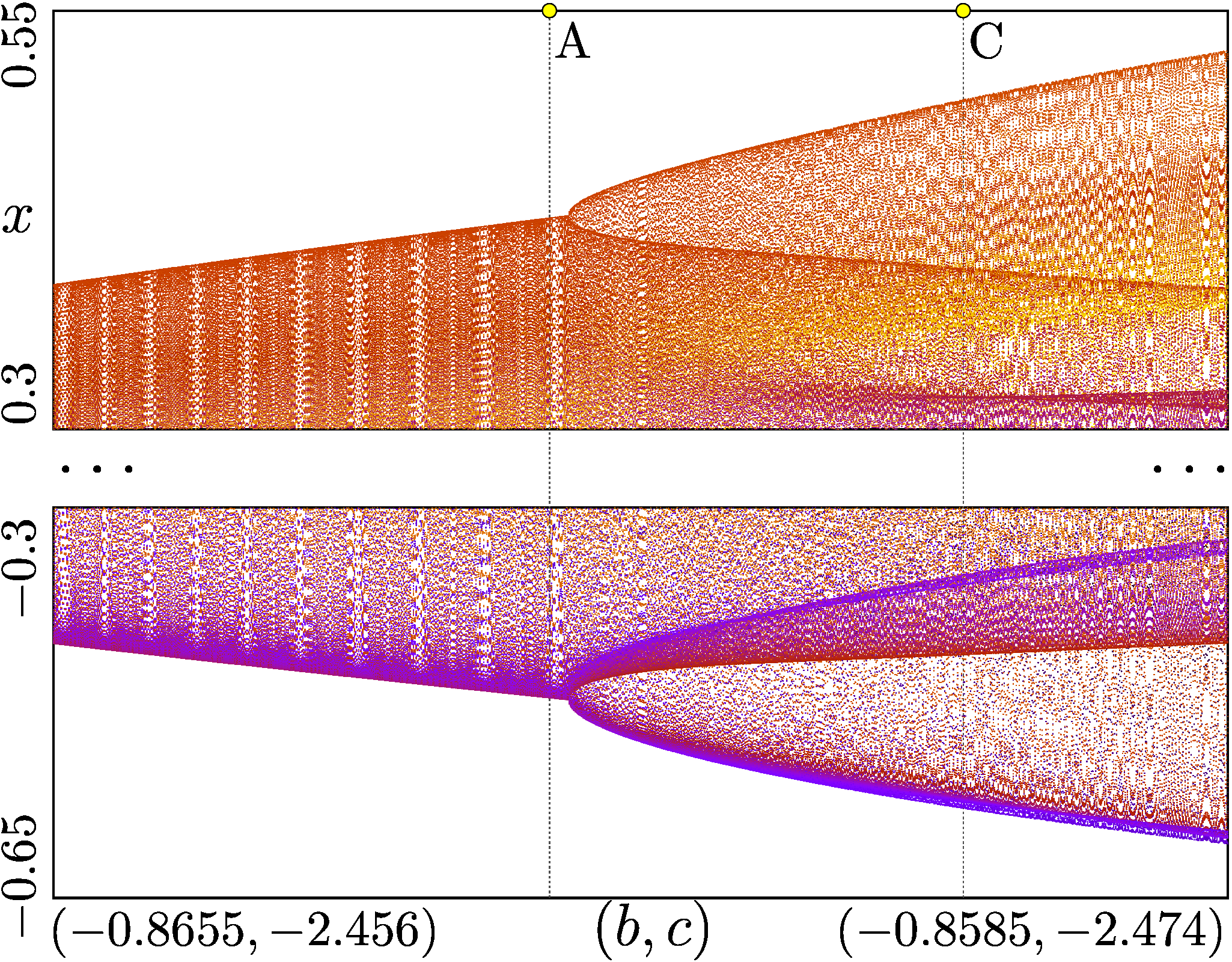}{d}
  \caption{(a) Parameter plane $(b,c)$ of the Mira map with $a=-2.5$, showing the overall
  bifurcation structure.  Gray indicates divergence, yellow indicates a fixed point, and green indicates no detected periodicity.  (b)
    Magnification of the rectangle marked in (a), showing the region
    near the 5-tongue. Parameter paths $\psi_1$ and $\psi_2$ exhibit the doubling of a resonant and a quasiperiodic invariant curves, as shown in (c) and (d) respectively. Parameter values
    corresponding to the configurations shown in 
    Figs.~\ref{fig:mira:resonant} and~\ref{fig:mira:quasiperiodic} are marked. 
    The points $A$, $C$, and $D$ in (c) refer to Figs.~\ref{fig:mira:resonant}(a), (c), and (d), respectively.
    The points $A$ and $C$ in (d) refer to ~\ref{fig:mira:quasiperiodic}(a) and (c), respectively.
  \label{fig:mira:2D}
  }
\end{figure*}

\begin{figure*}[tbh]
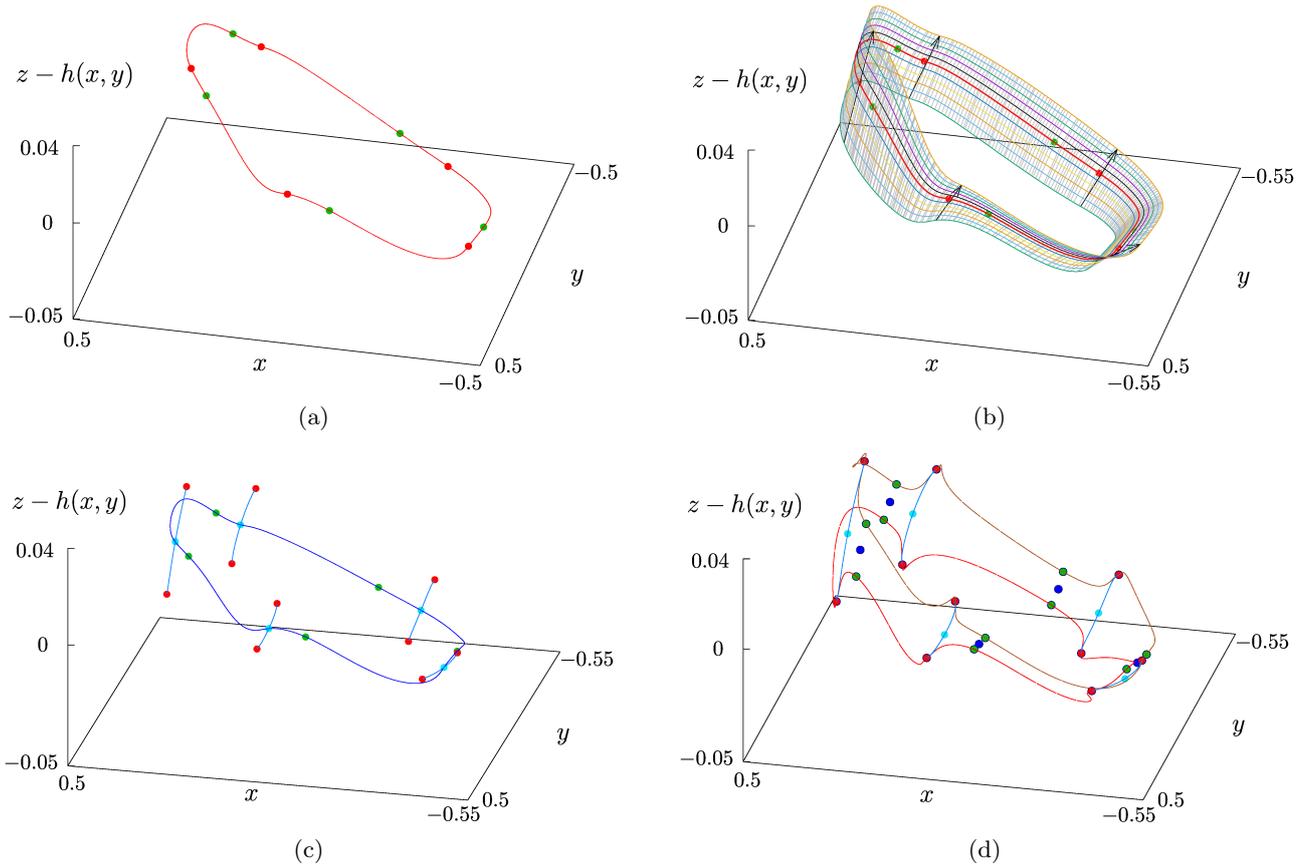
\centering
  \Includesubgraphics{.45\textwidth}{fig_Mira_r0Da}{a} \hspace{0.2in}
  \Includesubgraphics{.45\textwidth}{fig_Mira_r0Db}{b}\\ \vspace{0.1in}
  \Includesubgraphics{.45\textwidth}{fig_Mira_r0Dc}{c} \hspace{0.2in}
  \Includesubgraphics{.45\textwidth}{fig_Mira_r0Dd}{d}
  \caption{Loop-doubling of a resonant ICC in the
    Mira map.  (a) Attracting ICC before the doubling; (b) Doubling tangent space at the same parameter values as
    in (a); (c) Repelling ICC after the doubling of the node but
    before the doubling of the saddle;  (d) Attracting ICC consisting
    of two loops after both the node and the saddle have bifurcated. Additionally, in
    (c) and (d), unstable manifolds of the saddle $P^{uss}_5$ are also shown.
    Parameter values: $a=-2.5$; (a), (b) $b=-0.85578$, $c=-2.45869$;
    (c) $b=-0.8545207$, $c=-2.4613842$;
    (d) $b=-0.8533$, $c=-2.464$. The corresponding points are indicated in Fig.~\ref{fig:mira:2D}(c)
     by A, B, and C, respectively. To improve visibility, the ICCs and the doubling
tangent space are shown relative to
the surface $h(x,y)=a_1 x^2 + a_2 y^2 + a_3 xy + a_4 x + a_5 y + a_6$
with the coefficients $a_1$ through $a_6$ obtained by least-square
fitting of the ICC. 
    \label{fig:mira:resonant}
    }
\end{figure*}

\subsection{The Kamiyama maps}

Next, we consider two systems used in \cite{kamiyama2014classification}, namely:
\begin{equation}
\label{map:moebius}
   \begin{split}
    x_{n+1} &=R \cos {\theta} \cdot x_n 
        - R \sin {\theta} \cdot y_n \\ 
    &-(R-0.9) \cdot (x_n^2 + y_n^2) + E \cdot z_n\\
    y_{n+1} &= R \sin {\theta} \cdot x_n - R \cos {\theta}  \cdot y_n \\& -(R-0.9) \cdot (x_n^2 + y_n^2) \\
    z_{n+1} &= x_n       
   \end{split} 
\end{equation}
and
\begin{equation}
\label{map:cylinder}
   \begin{split}
    x_{n+1} &=R \cos {\theta} \cdot x_n 
        - R \sin {\theta} \cdot y_n \\ 
    &-(R-0.9) \cdot (x_n^2 + y_n^2) + E \cdot z_n\\
    y_{n+1} &= R \sin {\theta} \cdot x_n + R \cos {\theta}  \cdot y_n \\& -(R-0.9) \cdot (x_n^2 + y_n^2) \\
    z_{n+1} &= C \tanh z_n + D      
   \end{split} 
\end{equation}
where $\theta$ is an angle in degrees.

\section{Results}

\subsection{Loop doubling of a resonant torus}

Fig.~\ref{fig:mira:2D}(a) shows a 2D bifurcation diagram for the Mira map (\ref{eq:mira:map}) in the $(b,c)$ parameter plane with $a$ fixed at $-2.5$, which shows the bifurcation curve $\phi_{\text{\tiny NS}}$. Several periodicity
tongues originate from this bifurcation curve, including the period-5 tongue
magnified in Fig.~\ref{fig:mira:2D}(b). 
The figure also shows two parameter paths marked by $\psi_1$ (inside the resonant tongue) and $\psi_2$ (outside the tongue), corresponding to the doubling of a resonant ICC and a quasiperiodic ICC. The corresponding one-parameter bifurcation diagrams are
shown in Fig.~\ref{fig:mira:2D}(c) and~(d), respectively.

As can be seen in Fig.~\ref{fig:mira:2D}(c), in this case, the node doubles first. The stable period-5 cycle undergoes a flip bifurcation, leading to the appearance of a stable 10-cycle.
According to the Neimark-Sacker theorem, we
can expect that, sufficiently close to its origin, the considered period-5 
tongue contains a resonant ICC.

This is confirmed by Fig.~\ref{fig:mira:resonant}(a), which shows the structure of the phase space corresponding to the point marked A in Fig.~\ref{fig:mira:2D}(c). It shows an attracting ICC just before the doubling bifurcation: a union of the stable 5-node $P^{sss}_5$, the 5-saddle $Q^{uss}_5$, and its 1D unstable manifold.
In the notation used here, the three symbols in the superscript indicate the stability of the cycles's three multipliers, while the subscript refers to its period.  

To determine the type of doubling
bifurcation for this invariant curve, we compute its doubling tangent eigenspace by approximating it with the eigenvectors of the Jacobian $J_{f^5}$ corresponding to eigenvalues approaching $-1$. The eigenvectors have been obtained at a large number of points on the ICC and their tips have been joined by a green line, and different points on each arrow have been joined by lines of different color. In this way, the structure of the doubling tangent space becomes evident. As
illustrated in Fig.~\ref{fig:mira:resonant}(b), this tangent
space has a cylindrical topology.  This implies that a loop-doubling
bifurcation must occur.

The loop-doubling process is illustrated in
Figs.~\ref{fig:mira:resonant}(c) and (d). First, a flip bifurcation
transforms the stable 5-node $P^{sss}_5$ into a 5-saddle $P^{uss}_5$,
leading to the appearance of a stable 10-node, $P^{sss}_{10}$.  The
5-saddle $Q^{uss}_5$ remains unaffected by this bifurcation. As a
result, we obtain a repelling ICC composed of the 5-saddles
$P^{uss}_5$ and $Q^{uss}_5$, along with the 1D unstable manifold of
the latter (see Fig.~\ref{fig:mira:resonant}(c), which corresponds to the schematic illustration in Fig.~\ref{resonant-bif}(d)). The 1D unstable
 manifold of $P^{uss}_5$ approaches the stable 10-node $P^{sss}_{10}$ that is, at this stage, not yet involved in any ICC.

Next, the 5-saddle $Q^{uss}_5$ undergoes a flip bifurcation,
transforming into a 5-saddle $Q^{uus}_5$, while a 10-saddle
$Q^{uss}_{10}$ appears around it.  The unstable manifold of this
saddle converges to the stable 10-node $P^{sss}_{10}$, forming
an ICC.  As predicted, the resulting ICC consists of two loops invariant under the second iterate $f^2$. The sequence of events follows the route shown in Fig.~\ref{resonant-bif} (a) $\rightarrow$ (d) $\rightarrow$ (e) $\rightarrow$ (f).

The eigenvalues of $J_{f^5}$ corresponding to the parameter values of Fig.~\ref{fig:mira:resonant}(a) corresponding to the doubling, tangent and third directions are $-0.99767964$, $0.86129688$, and $0.53415424$ respectively. The third eigenvalue is positive, which conforms to the Gardini-Sushko conjecture mentioned in Section~\ref{sec:gardini}.

\begin{figure}[tbh]
\centering
\Includesubgraphics{0.4\textwidth}{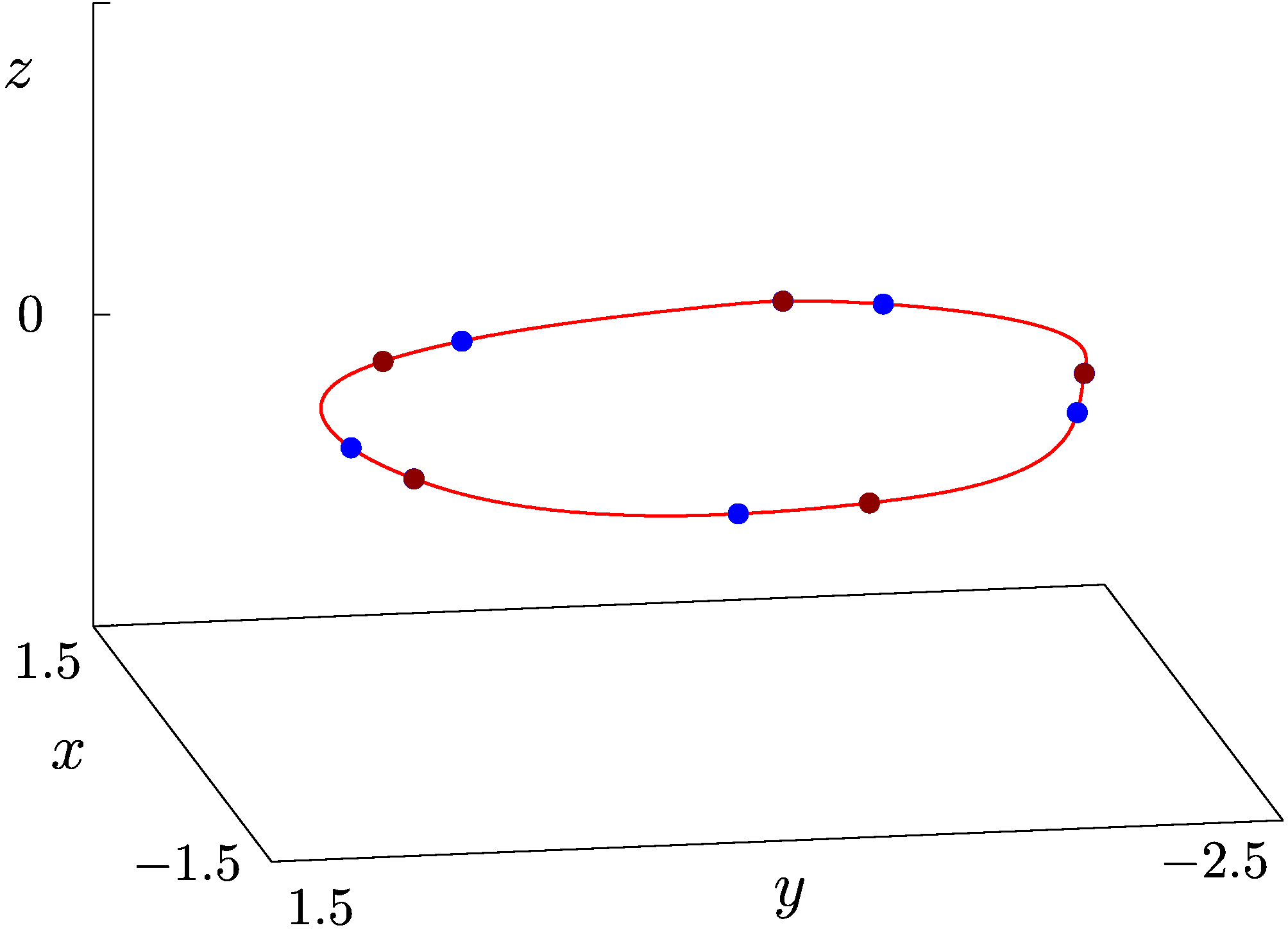}{a}\\
\Includesubgraphics{0.4\textwidth}{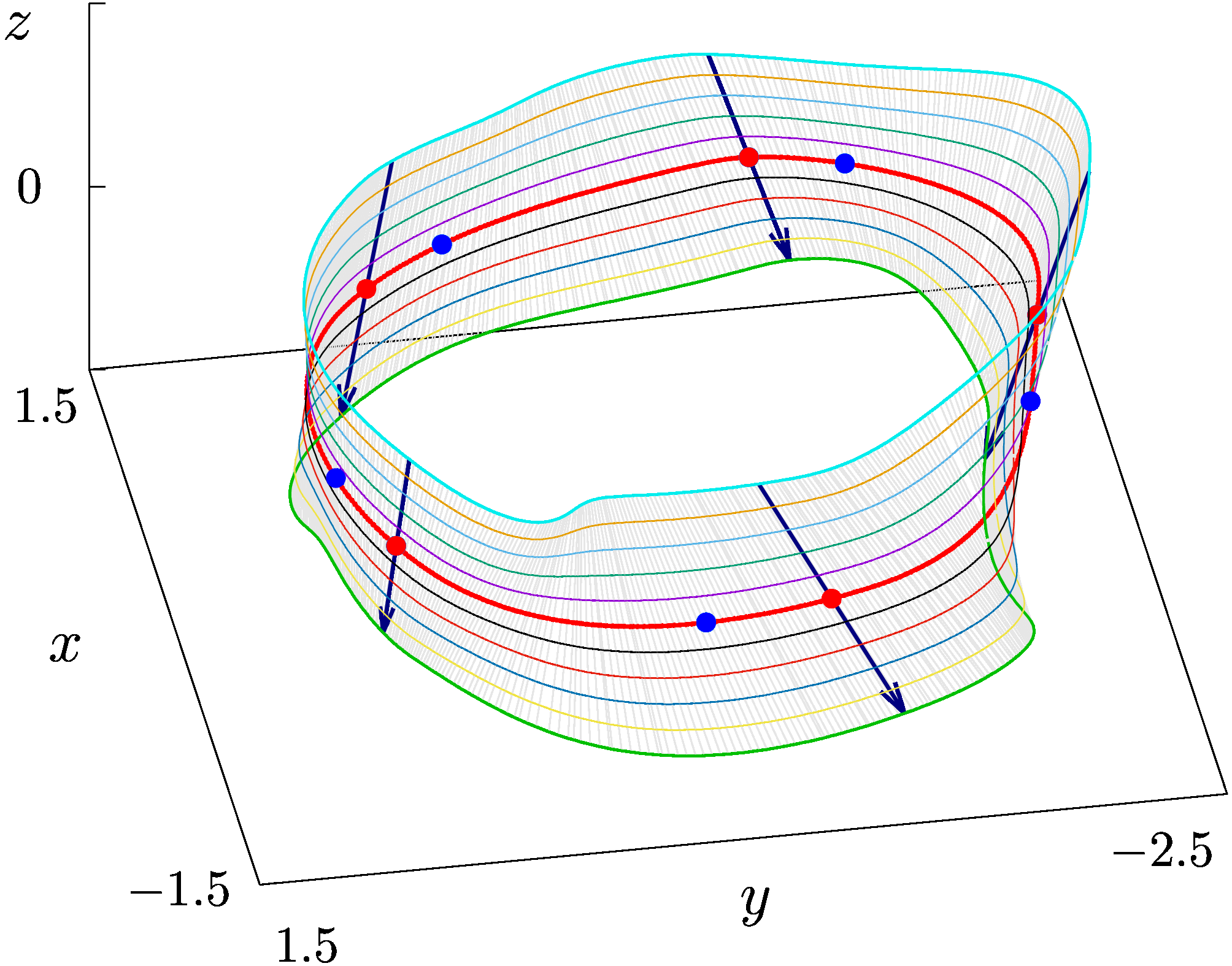}{b}\\
\Includesubgraphics{0.4\textwidth}{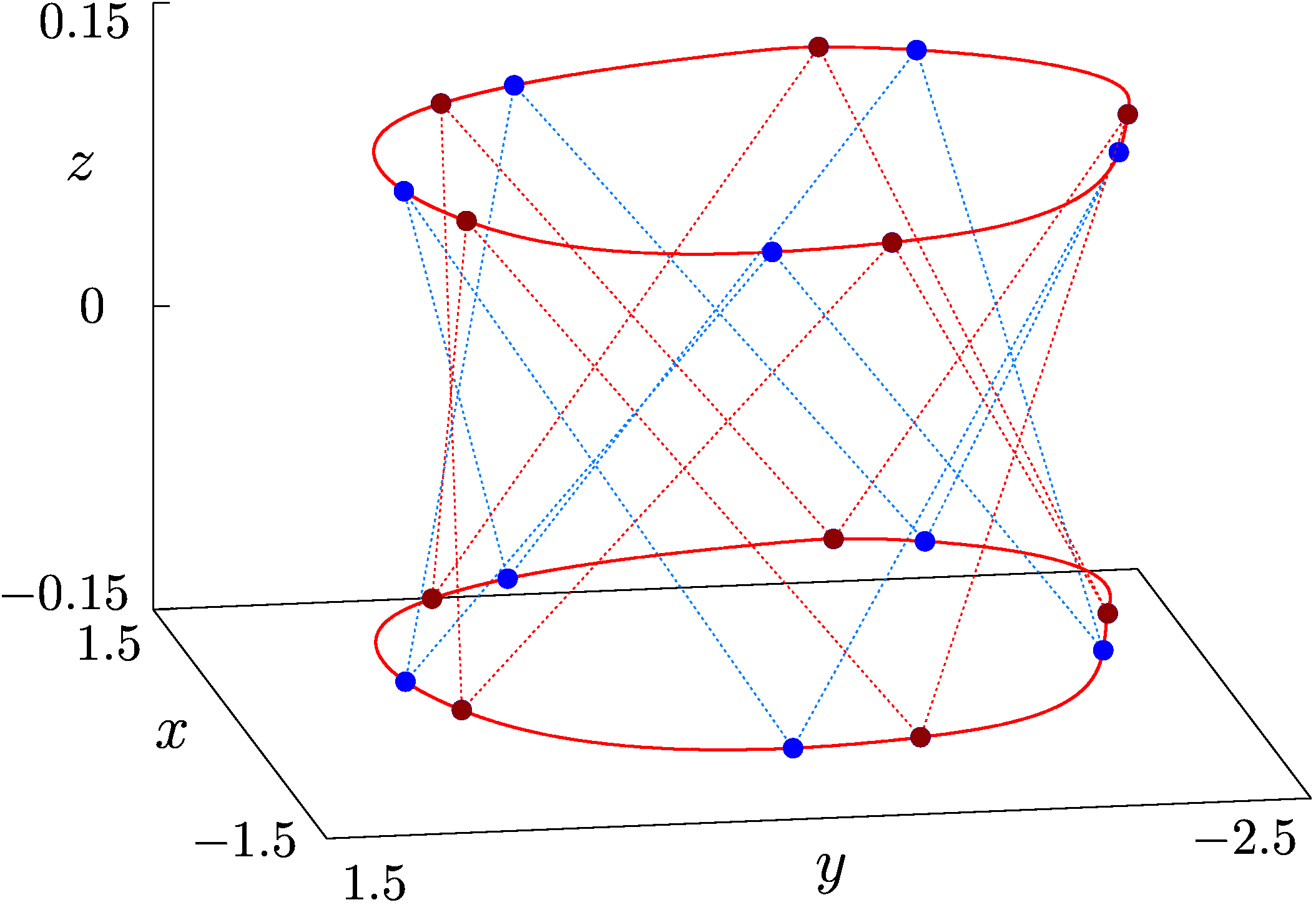}{c}

\caption{\label{fig:loop:doubling}Loop doubling in map~(\ref{map:cylinder}), (a) before the bifurcation; (b) The tangent eigenspace showing the eigenvectors of $J_{f^5}$ corresponding to the doubling eigenvalue evaluated at the nodes; (c) the ICC after the bifurcation. The red dots represent a period-5 stable cycle and the blue dots represent a period-5 saddle. The saddle-node connection via the unstable manifolds of the saddle forms the ICC, shown in red.  Parameters: $R=1.15$, $E=0.1$, $\theta=82$, $D=0$, (a) and (b) $C=-0.999$, (c) $C=-1.03$.
}
\end{figure}

However, for the Kamiyama map (\ref{map:cylinder}), we get a different result. For this system, the loop doubling phenomenon is shown in Fig.~\ref{fig:loop:doubling}. In this case also, the tangent eigenspace (Fig.~\ref{fig:loop:doubling}(c)) shows a cylindrical structure.
For the ICC in Fig.~\ref{fig:loop:doubling}(a), the eigenvalues of $J_{f^5}$ corresponding to the doubling, tangent and third directions evaluated at the points of the stable cycle are
$-0.99501, 0.72672,-0.06028$, respectively. Notice that the third eigenvalue is negative, which violates the Gardini-Sushko conjecture.

\subsection{Loop doubling of an ergodic torus}

\begin{figure*}[t]\centering
  \Includesubgraphics{.48\textwidth}{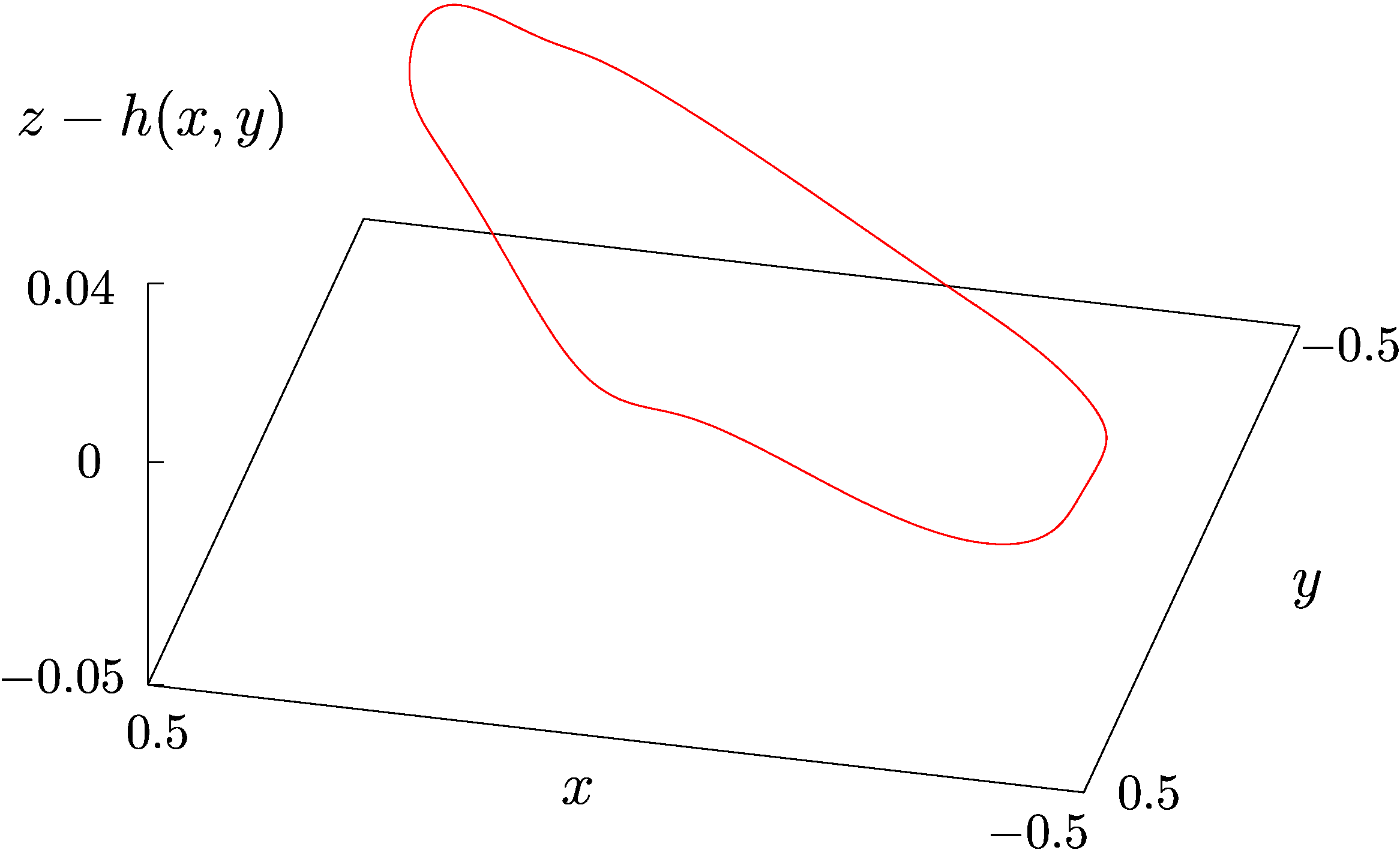}{a} \\
  \Includesubgraphics{.48\textwidth}{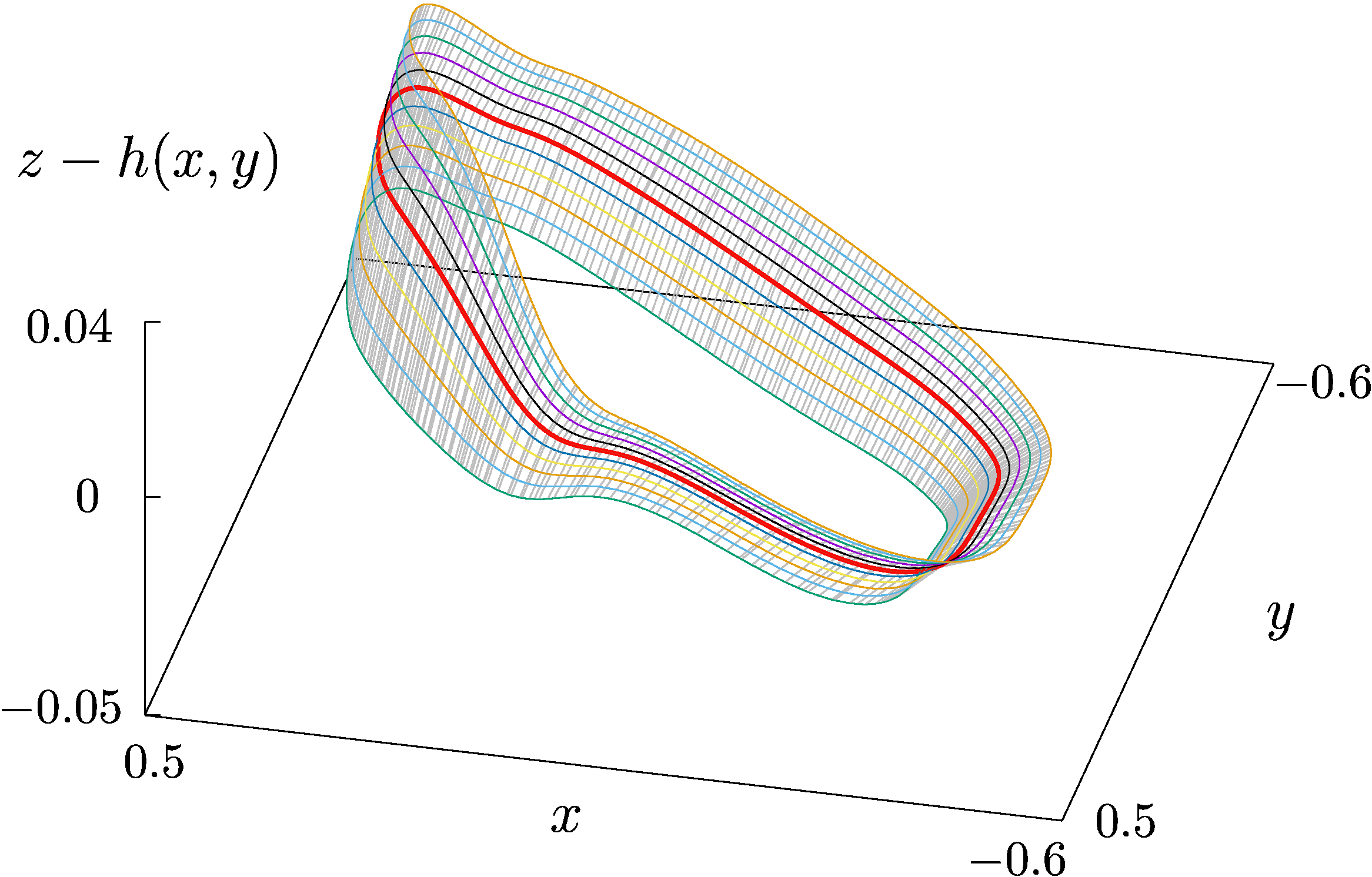}{b}\\
  \Includesubgraphics{.48\textwidth}{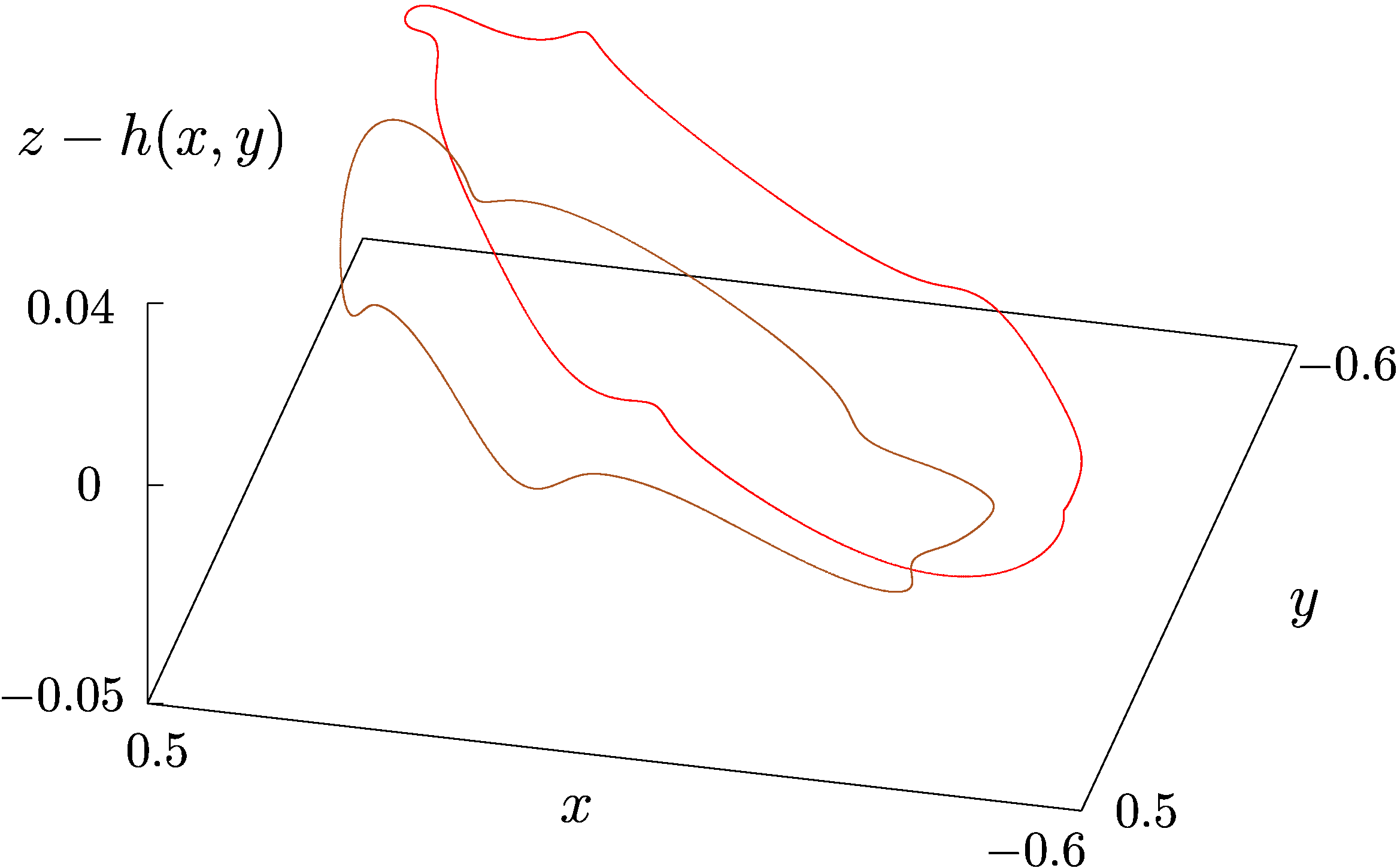}{c}
  \caption{Loop-doubling of a sufficiently quasiperiodic invariant closed curve in
    the Mira map.  (a) Attracting ICC before the doubling; (b) Doubling tangent space at the same parameter values as in (a); (c) Attracting ICC consisting of two loops after the doubling.     Parameter values: $a=-2.5$; (a), (b) $b=-0.862501$, $c=-2.463746$; (c) $b=-0.86$, $c=-2.47$; The corresponding points are indicated in Fig.~\ref{fig:mira:2D}(d) as A and C, respectively.  }
  \label{fig:mira:quasiperiodic}
\end{figure*}

To illustrate the doubling of a quasiperiodic ICC\footnote{One should note that a finite-precision computer cannot identify a truly irrational frequency ratio. So, when we say an orbit is quasiperiodic, one should take it in the sense of being `sufficiently quasiperiodic'.}, consider the
bifurcation diagram for the Mira map, shown in Fig.~\ref{fig:mira:2D}(d), which
corresponds to the parameter path marked in Fig.~\ref{fig:mira:2D}(b)
as $\psi_2$.  In Fig.~\ref{fig:mira:2D}(d), only the upper and
lower parts of the bifurcation diagram are shown (with the middle part
omitted), making the doubling clearly visible.  The sufficiently
quasiperiodic ICC at the parameter point marked as A in
Fig.~\ref{fig:mira:2D}(d) immediately before the doubling bifurcation
is shown in Fig.~\ref{fig:mira:quasiperiodic}(a).  To determine the
type of doubling the ICC undergoes, we need to choose the iterate
$f^p$ to compute the tangent space along the ICC.    To this end, we numerically estimate the
rotation number $\omega$ of the ICC, obtaining $\omega\approx
0.39860585$. Then, we consider $\frac{2}{5}$ as a sufficiently close
approximation of $\omega$ and calculate the tangent eigenspace in the
doubling direction using the fifth iterate $f^5$. Note that the same
procedure, namely a rational approximation of the numerically computed
rotation number, can also be used to compute the eigenvalues in the second Poincar\'e section discussed in~\cite{banerjee2012local}.

The doubling tangent space, approximated using the doubling eigenvectors of
$J_{f^5}$ at points along the quasiperiodic ICC, is shown
in Fig.~\ref{fig:mira:quasiperiodic}(b). As clearly visible, it
has a cylindrical structure. Therefore, the doubling bifurcation that
the considered ICC undergoes is a loop
doubling. This is confirmed by Fig.~\ref{fig:mira:quasiperiodic}(c),
which shows the ICC after the doubling.

\subsection{Length doubling of a resonant torus}

\begin{figure*}[t!]
\centering
  \Includesubgraphics{.4\textwidth}{fig_Kamiyama_2Da}{a}\hspace{0.3in}
  \Includesubgraphics{.4\textwidth}{fig_Kamiyama_2Db}{b}\\

  \medskip
    \Includesubgraphics{.4\textwidth}{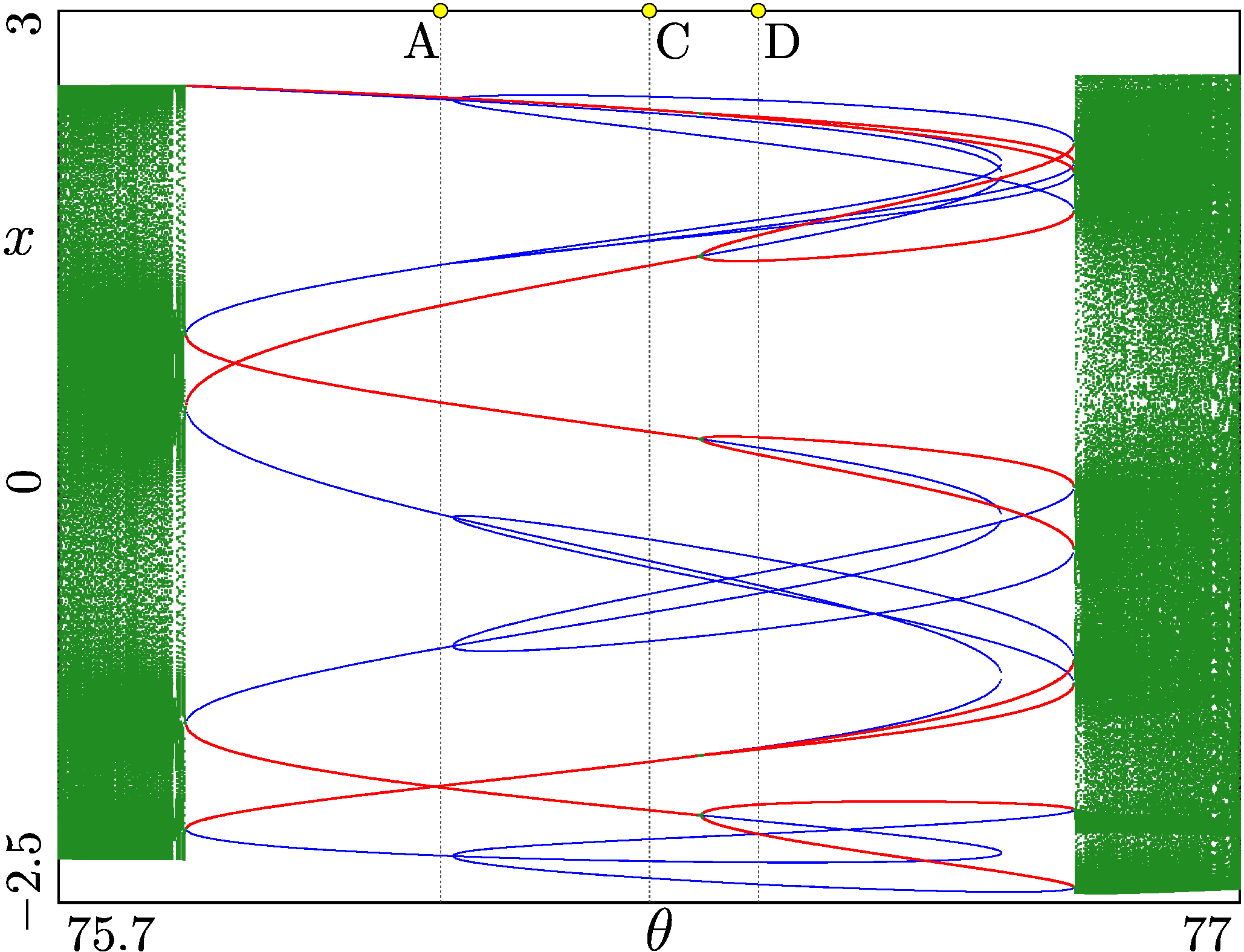}{c}\hspace{0.3in}
  \Includesubgraphics{.4\textwidth}{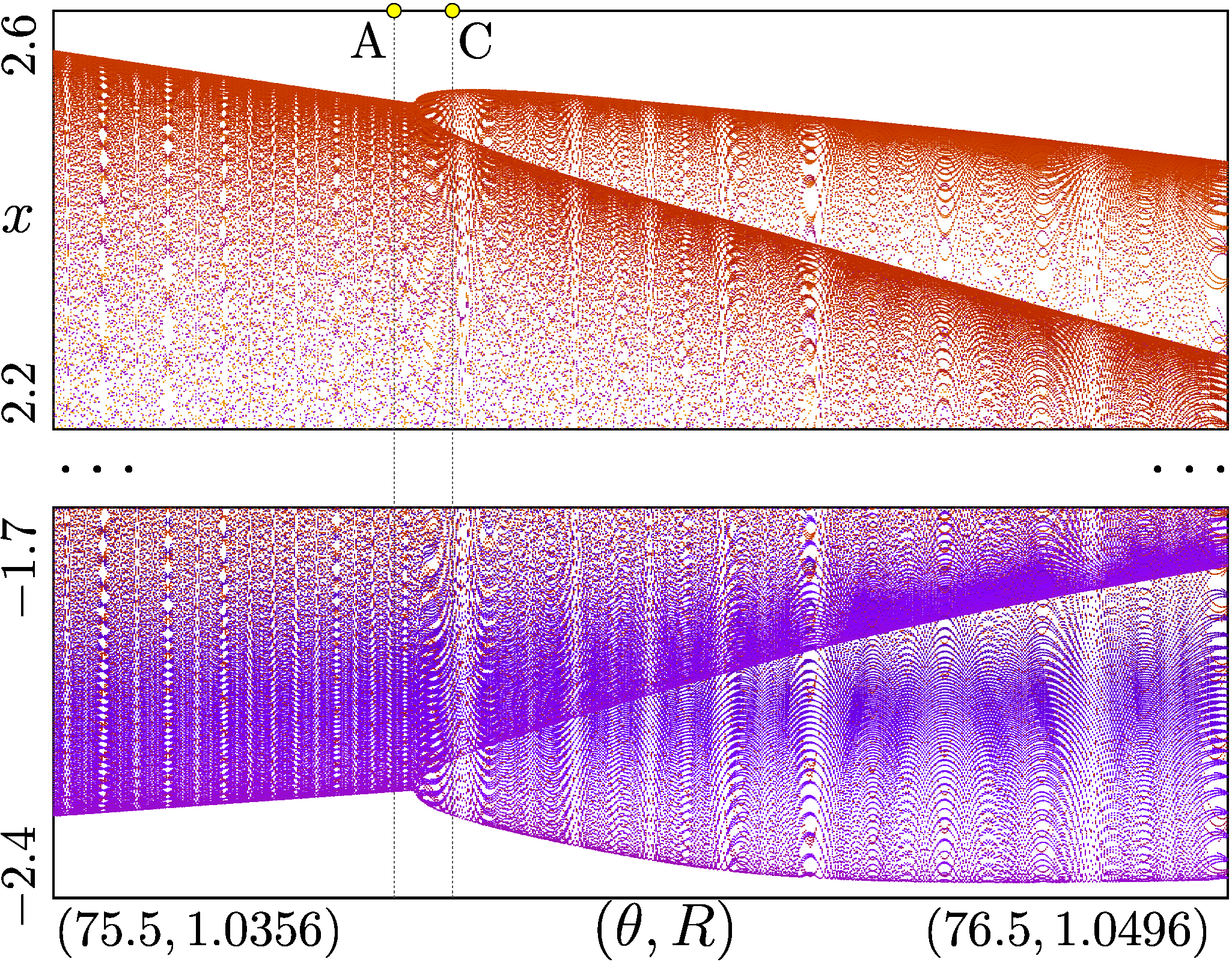}{d}
  \caption{(a) Parameter plane $(b,c)$ of the Kamiyama map (\ref{map:moebius}) with $E=-0.3$, showing the overall bifurcation structure.  Gray indicates divergence, blue indicates period~5 and green indicates period~10.  (b) Magnification of the rectangle marked in (a), showing the region near the 5-tongue. Parameter paths $\psi_1$ and $\psi_2$ exhibit the doubling of a resonant and a quasiperiodic invariant curves, as shown in (c) and (d), respectively. Parameter values corresponding to the configurations shown in 
    Figs.~\ref{fig:length:doubling} and~\ref{fig:length:doubling:q} are marked. The points $A$, $C$, and $D$ in (c) refer to Figs.~\ref{fig:length:doubling} (a), (c), and (d), respectively.
    The points $A$ and $C$ in (d) refer to Fig.~\ref{fig:length:doubling:q}(a) and (c), respectively.}
  \label{fig:kamiyama:2D} 
  
\end{figure*}

\begin{figure*}[t]
\centering
\Includesubgraphics{0.4\textwidth}{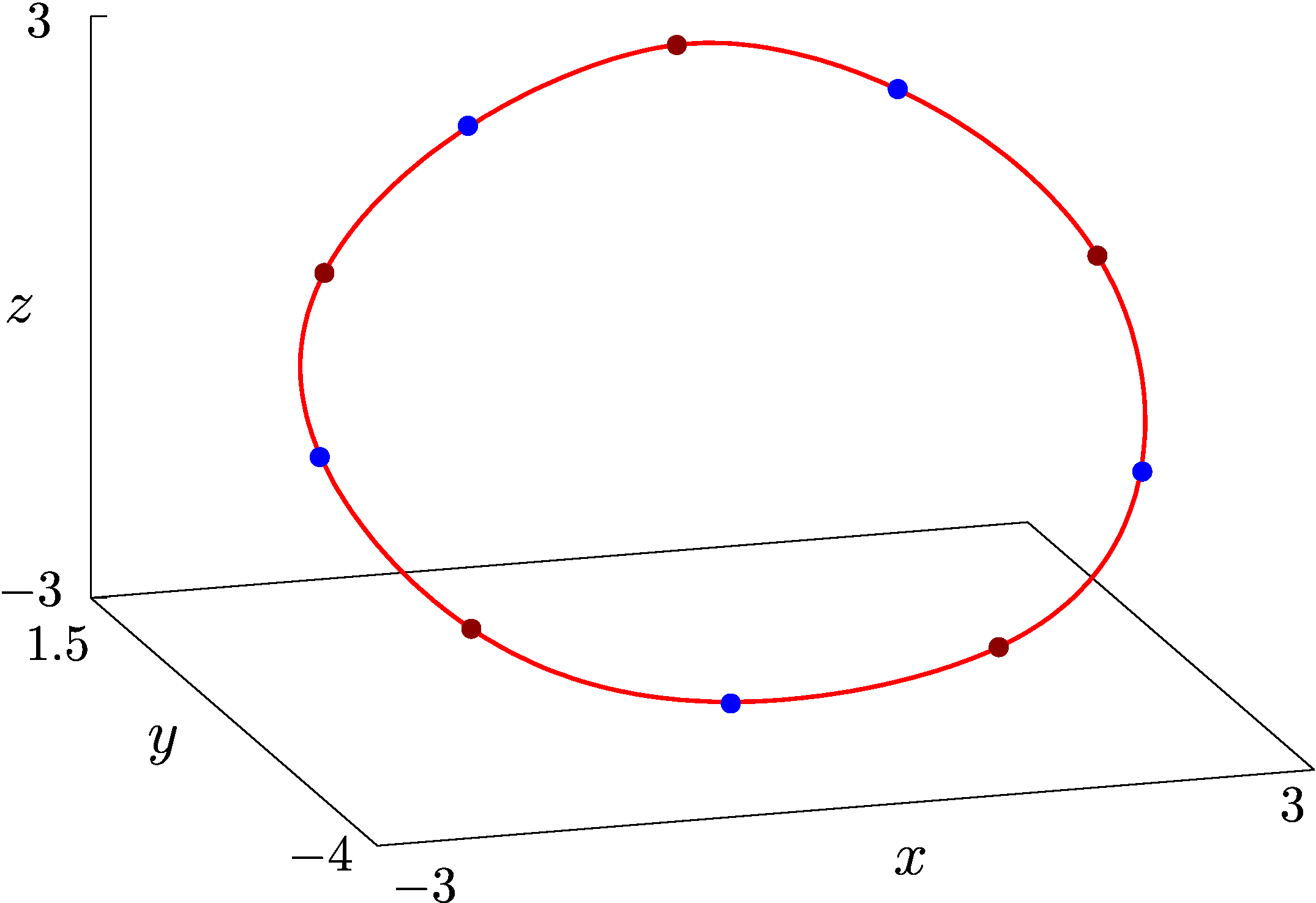}{a}\hspace{0.5in}
\Includesubgraphics{0.4\textwidth}{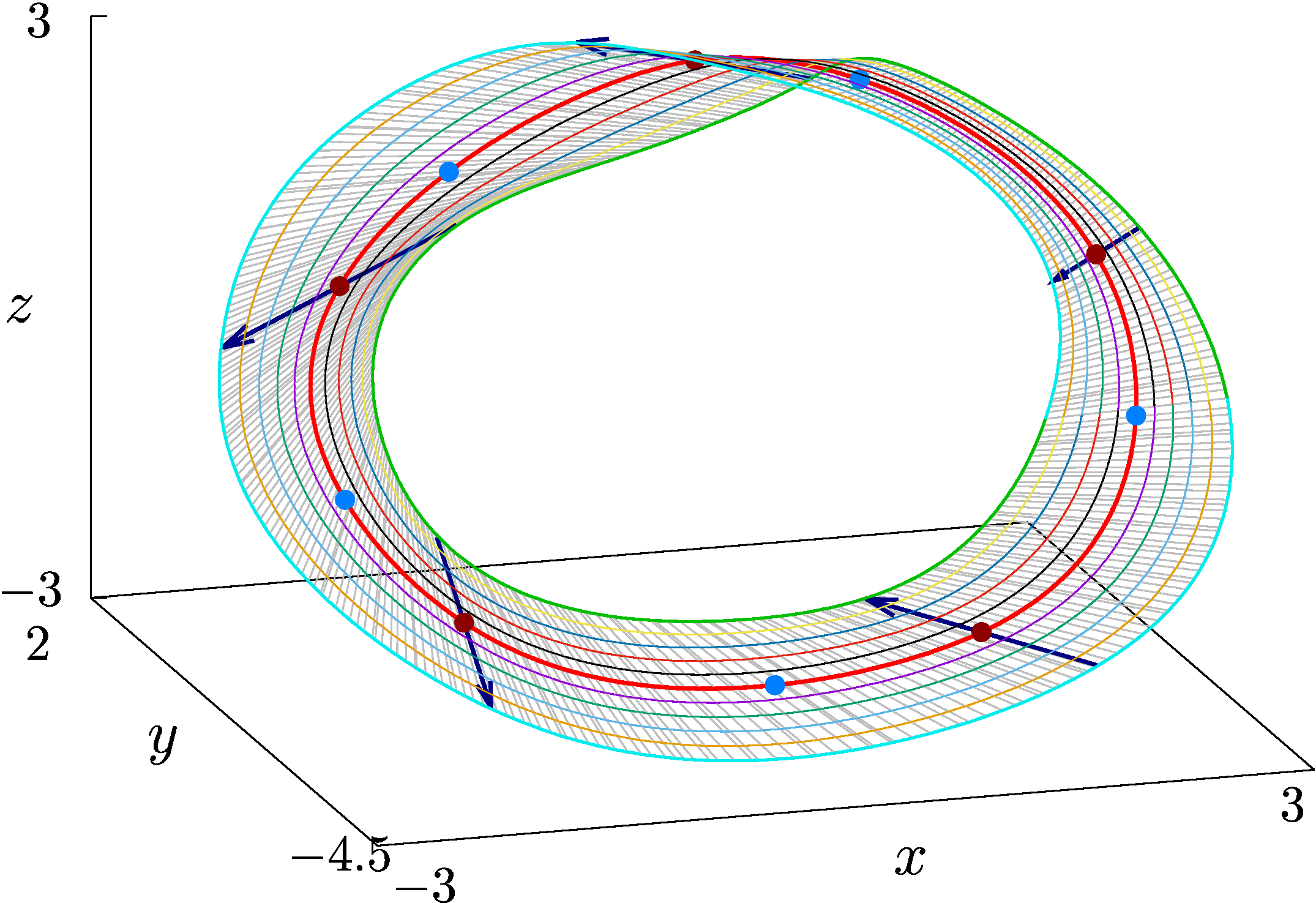}{b}\\
\Includesubgraphics{0.4\textwidth}{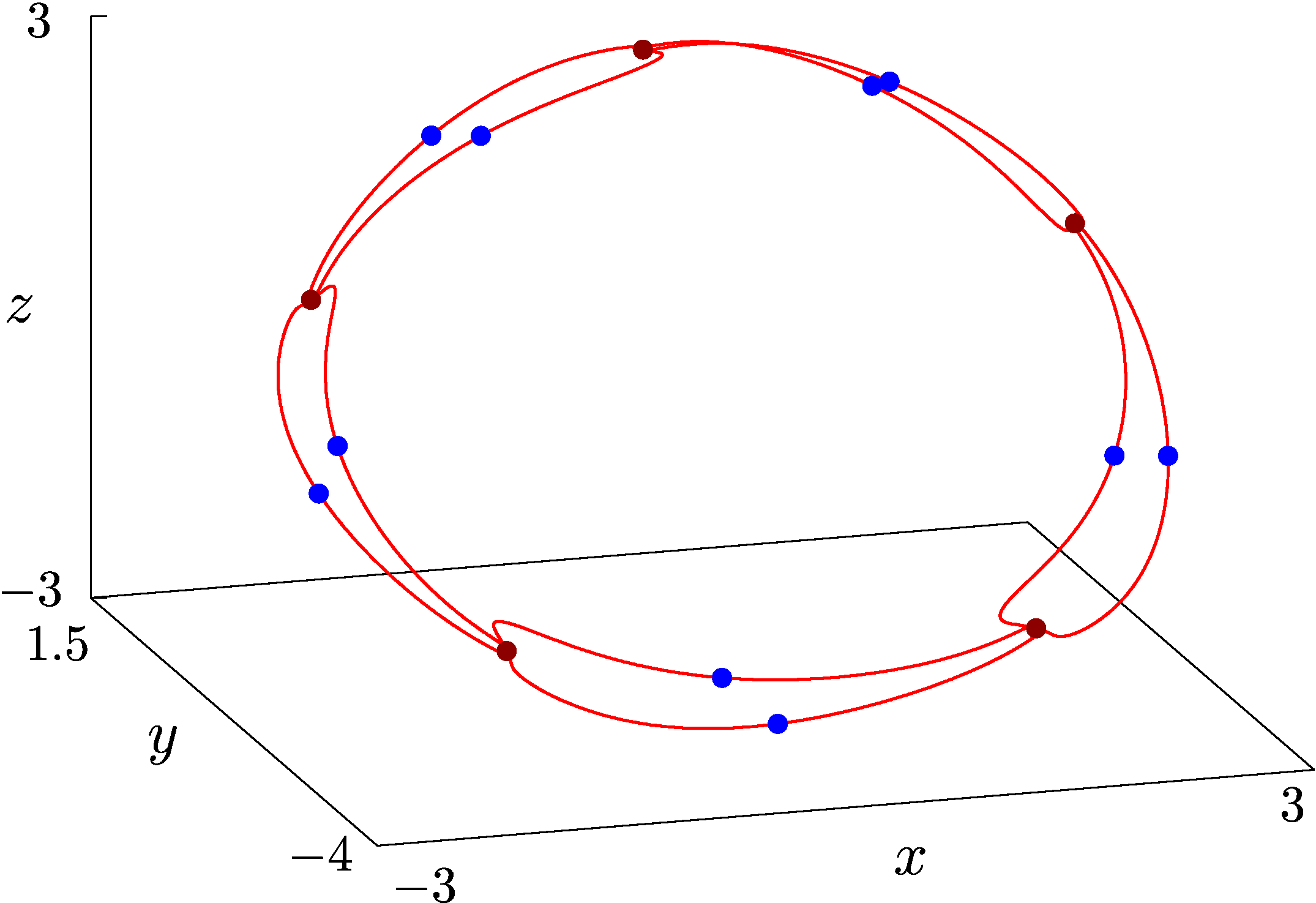}{c} \hspace{0.5in}
\Includesubgraphics{0.4\textwidth}{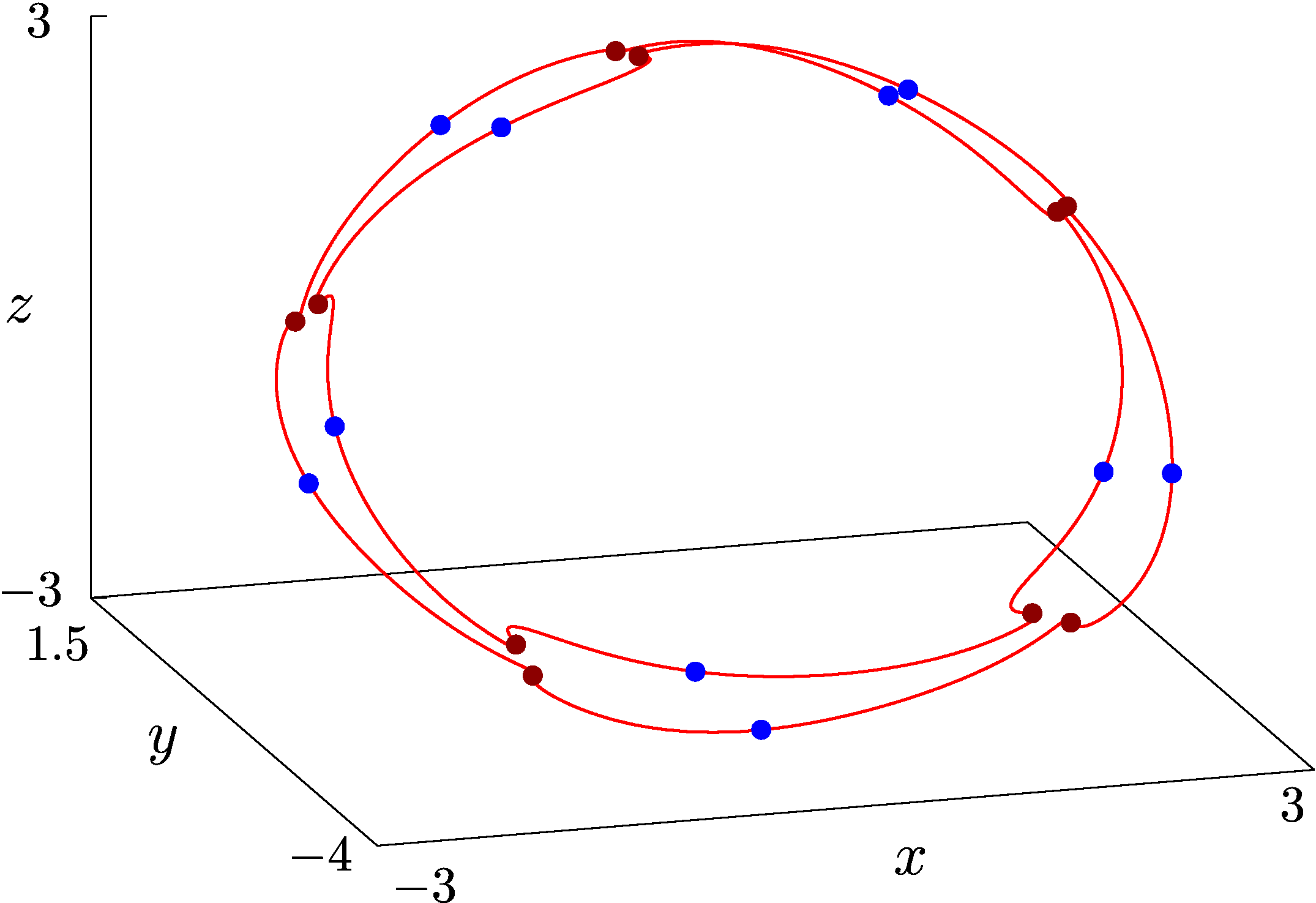}{d}
\caption{\label{fig:length:doubling}Length doubling in 
map~(\ref{map:moebius}). (a) The ICC before the bifurcation, (b) the tangent space obtained by the eigenvectors of $J_{f^5}$ corresponding to the doubling  
eigenvalue evaluated at the points of the ICC. The eigenvectors at the nodes are also shown. 
(c) Bi-layered ICC after doubling 
 bifurcation of the saddle but before
 doubling bifurcation of the node.
 (d) Period-10 ICC after doubling 
 bifurcation of the node. Figures (a), (c), and (d) refer to parameters marked as $A$, $B$, and $C$ in Fig.~\ref {fig:kamiyama:2D}(c), respectively.
Parameters: $R=1.038$, $E=-0.3$,
(a), (b) $\theta=76.12$,
(c) $\theta=76.35$,
(d) $\theta=76.47$.}
\end{figure*}

For the map (\ref{map:moebius}), Fig.~\ref{fig:kamiyama:2D}(a) shows a 2D bifurcation diagram showing several periodicity tongues, and a zoom on to the period-5 tongue is shown in Fig.~\ref{fig:kamiyama:2D}(b). It shows two paths $\psi_1$ and $\psi_2$ along which parameters are varied. 

Fig.~\ref{fig:length:doubling} shows the invariant closed curves before and after a length-doubling bifurcation as the parameter is varied along $\psi_1$. Fig.~\ref{fig:length:doubling}(a) shows the ICC before the bifurcation begins, corresponding to the parameter marked A in Fig.~\ref{fig:kamiyama:2D}(c). 
Applying our proposed method, we obtain the set of eigenvectors along the ICC. 
The tangent eigenspace related to the doubling eigenvalue before the bifurcation is shown in Fig.~\ref{fig:length:doubling}(b). 
It is clear that this tangent eigenspace has a M\"obius strip structure. According to our theory, there should be a length-doubling bifurcation as the parameter is varied.

Fig.~\ref{fig:kamiyama:2D}(c) shows that, in this case, the saddle doubles first. The situation after the saddle doubles but before the node doubles is shown in Fig.~\ref{fig:length:doubling}(c), for the parameter value marked B in Fig.~\ref{fig:kamiyama:2D}(c). 
 In the next stage, the node also undergoes a period-doubling bifurcation, resulting in the length-doubled ICC (Fig.~\ref{fig:length:doubling}(d)) for the parameter value marked C in Fig.~\ref{fig:kamiyama:2D}(c). Therefore, this ICC doubling occurs following the sequence illustrated as Fig.~\ref{resonant-bif}(a) $\rightarrow$ (b) $\rightarrow$ (c) $\rightarrow$ (f).

The eigenvalues of $J_{f^5}$ 
corresponding to the doubling, tangent, and third directions just before the bifurcation, evaluated at the points
of the stable cycle are $-0.95979, 0.85489, 3.58612\times 10^{-6}$, respectively.  We notice that the third eigenvalue is positive, which contradicts the Gardini-Sushko conjecture.

\subsection{Length Doubling of a quasiperiodic torus}

\begin{figure}[t]
\centering
\Includesubgraphics{0.4\textwidth}{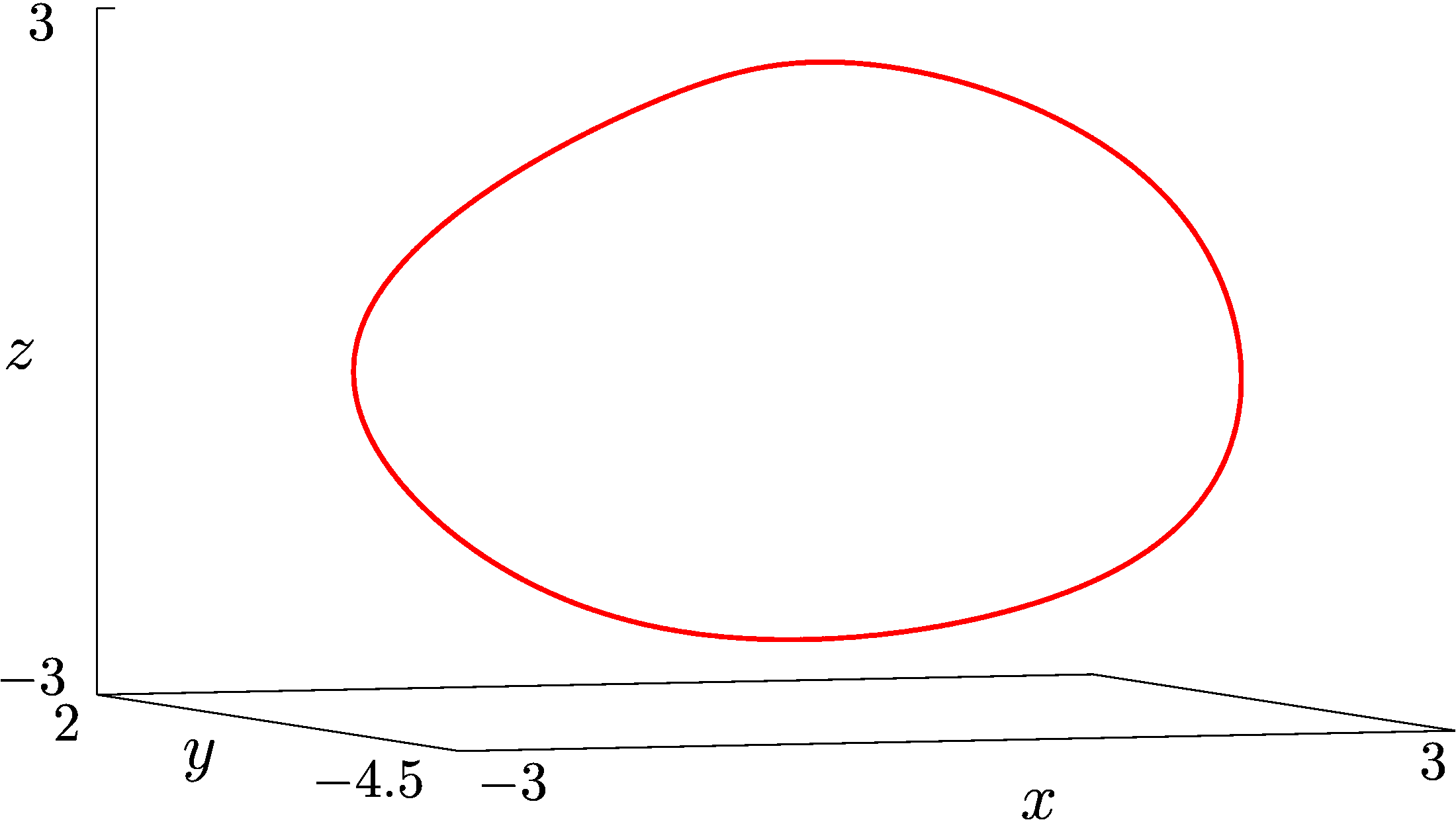}{a}\\
\Includesubgraphics{0.4\textwidth}{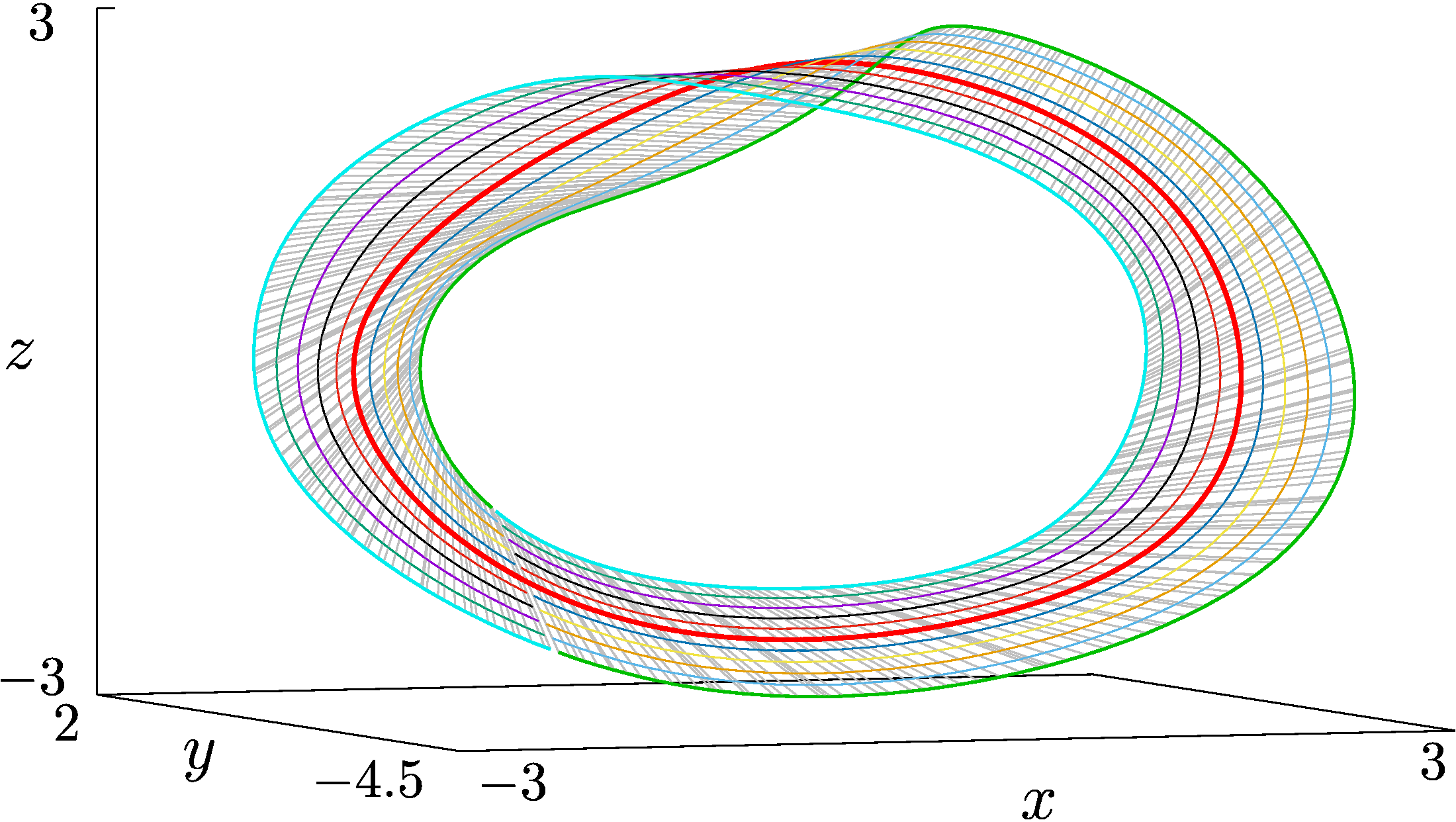}{b}\\
\Includesubgraphics{0.4\textwidth}{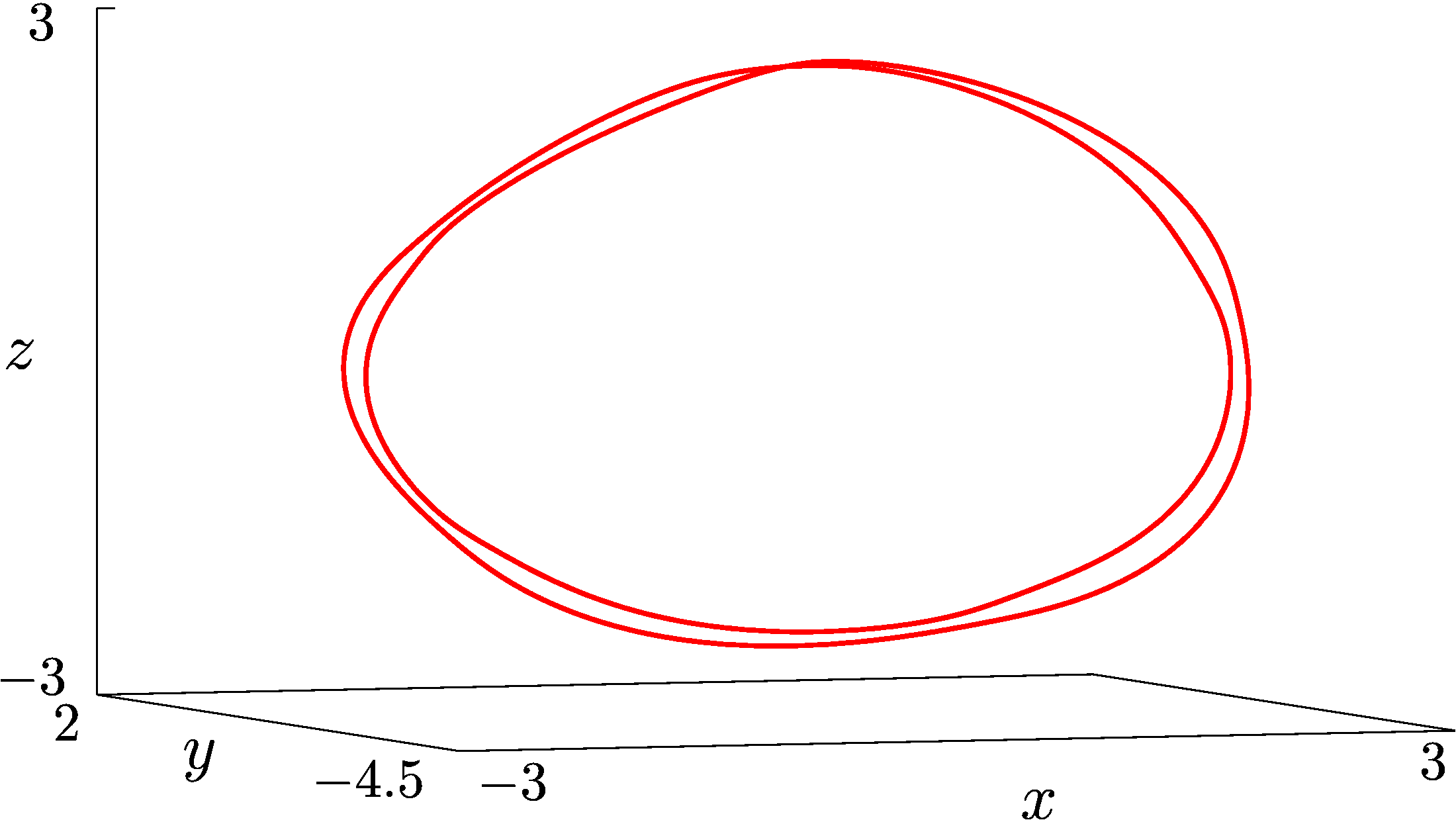}{c}
\caption{\label{fig:length:doubling:q}Length doubling of a sufficiently quasiperiodic ICC in map~(\ref{map:moebius}). (a) The ICC before the bifurcation, (b) the tangent eigenspace, and (c) the orbit after the bifurcation. The parameter values for (a) and (b) are $E=-0.3$, $R=1.0397$, $\theta=75.79$ (marked as A in Fig.~\ref{fig:kamiyama:2D}d) and those for (c) are $E=-0.3$, $R=1.0404$
$\theta=75.8398$ (marked as C in Fig.~\ref{fig:kamiyama:2D}d). Rotation number $\approx 0.199220$.
}
\end{figure}

Fig.~\ref{fig:length:doubling:q} shows the appearance of a double-length quasiperiodic ICC for the map (\ref{map:moebius}) as the parameter is varied along the line $\psi_2$ shown in Fig.~\ref{fig:kamiyama:2D}(b). 
Fig.~\ref{fig:length:doubling:q}(a) shows the quasiperiodic ICC before bifurcation. The tangent eigenspace (Fig.~\ref{fig:length:doubling:q}b) computed following the procedure outlined above shows a M\"obius strip structure. Finally, Fig.~\ref{fig:length:doubling:q}(c) shows the resulting length-doubled ICC.

\subsection{ICC with saddle-focus connection}

\begin{figure}[t]
\centering
\includegraphics[width=0.35\textwidth, height=0.25\textwidth]{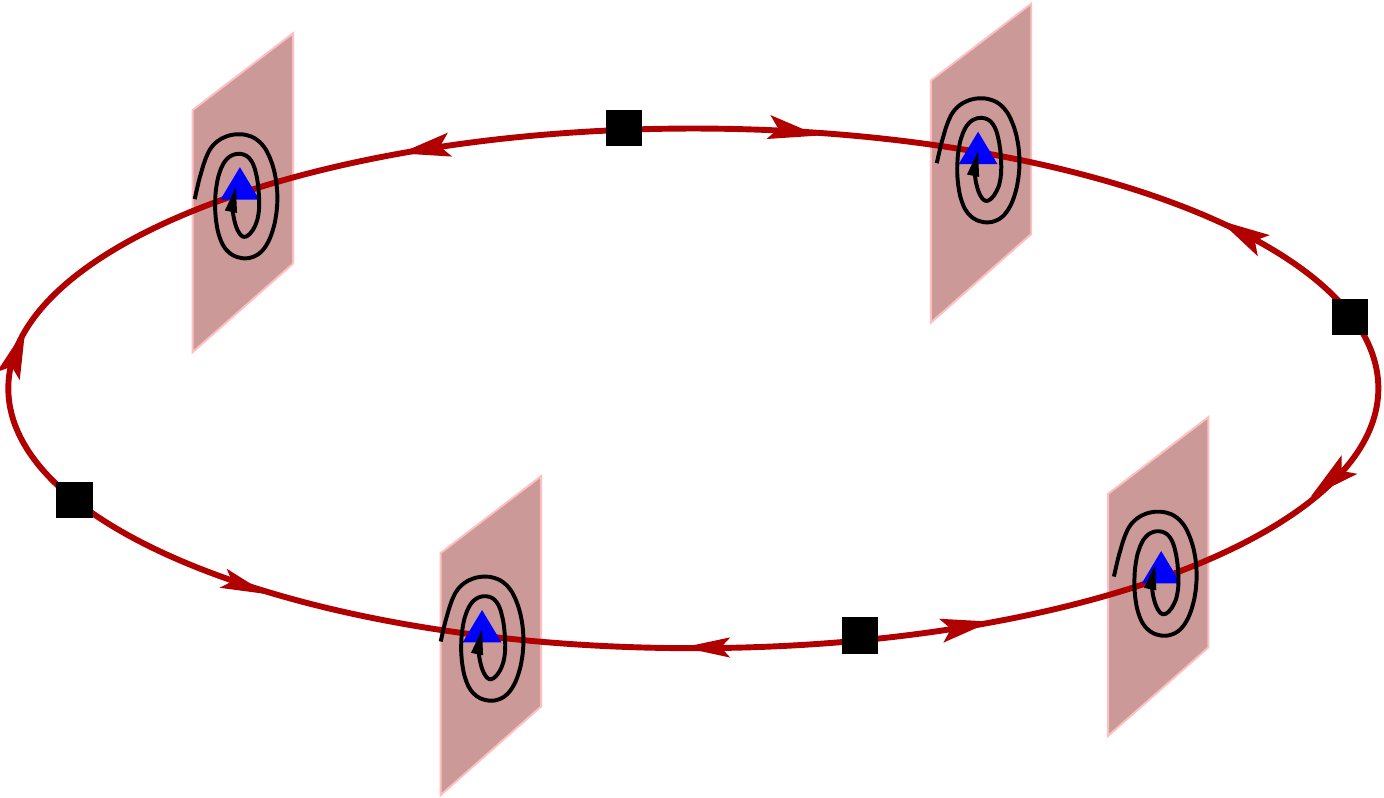}{\footnotesize (a)}
\includegraphics[width=0.35\textwidth, height=0.25\textwidth]{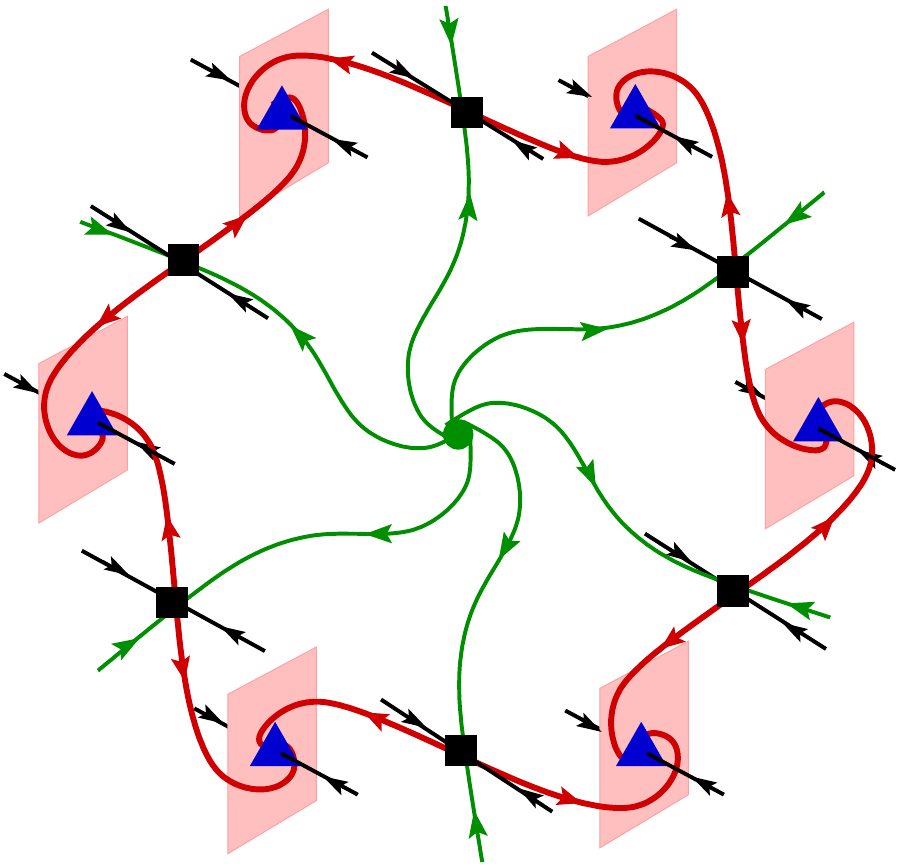}{\footnotesize (b)}

\caption{Schematic diagram showing two types of invariant closed curves connecting saddle cycles and focus-type nodes \label{saddle-focus}}
\end{figure}

The phenomena discussed in this paper occur in maps of dimension 3 or higher. If we consider a 3D system, there are three eigenvalues of the Jacobian matrix. All three eigenvalues can be real, or one can be real and the other two complex conjugate. In the latter case, the ICC could be formed in two possible ways:
\begin{enumerate}
\item The ICC is formed by the unstable manifolds of the saddle cycle that are tangent to the real eigenvectors;
\item The unstable manifolds of the saddles spiral onto the stable nodes (the saddle-focus connection) thus forming an ICC of infinite length. 
\end{enumerate}
Fig.~\ref{saddle-focus} shows the two types of structures schematically.

\begin{figure*}[t]
\centering
\Includesubgraphics{0.3\textwidth}{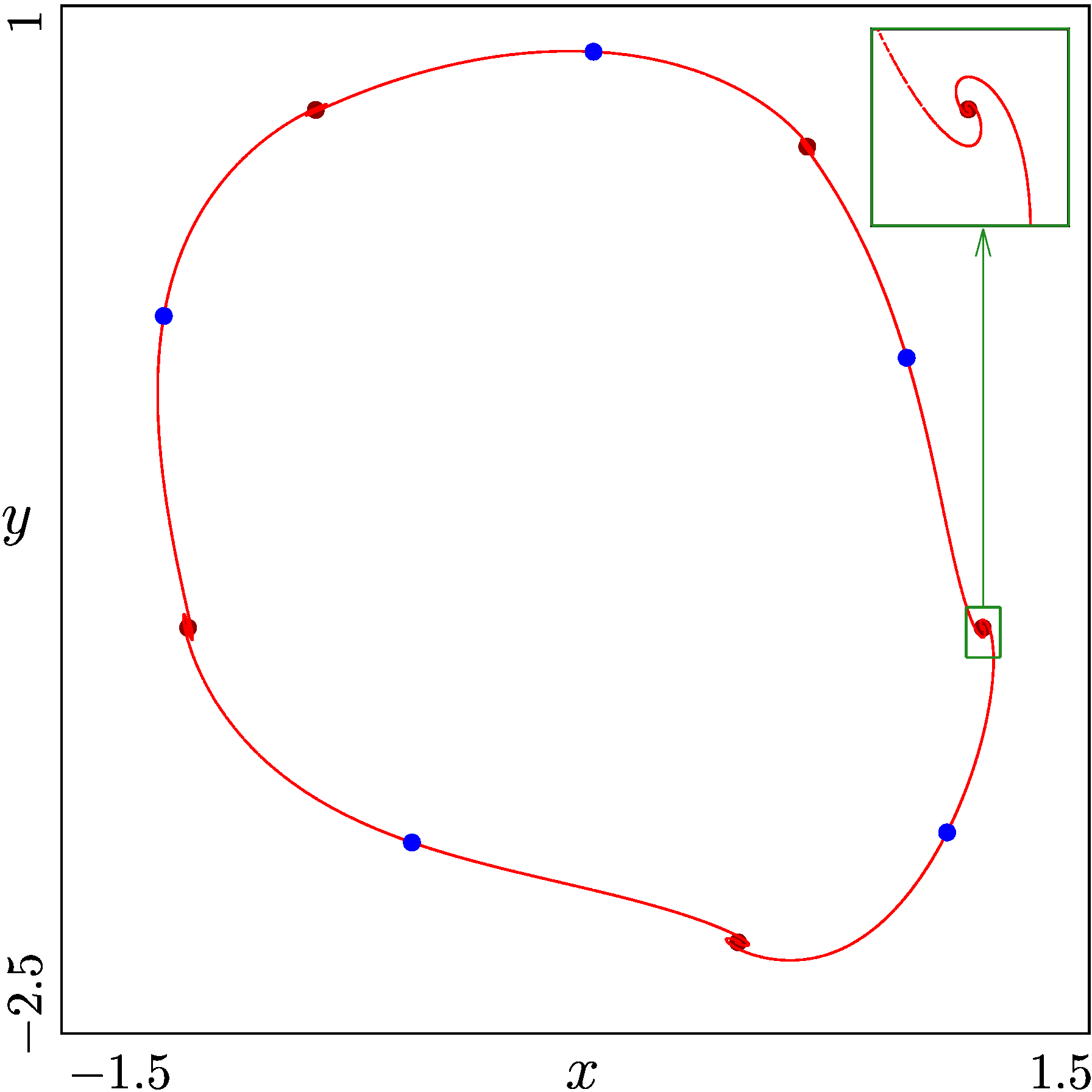}{a}\\ 
\Includesubgraphics{0.4\textwidth}{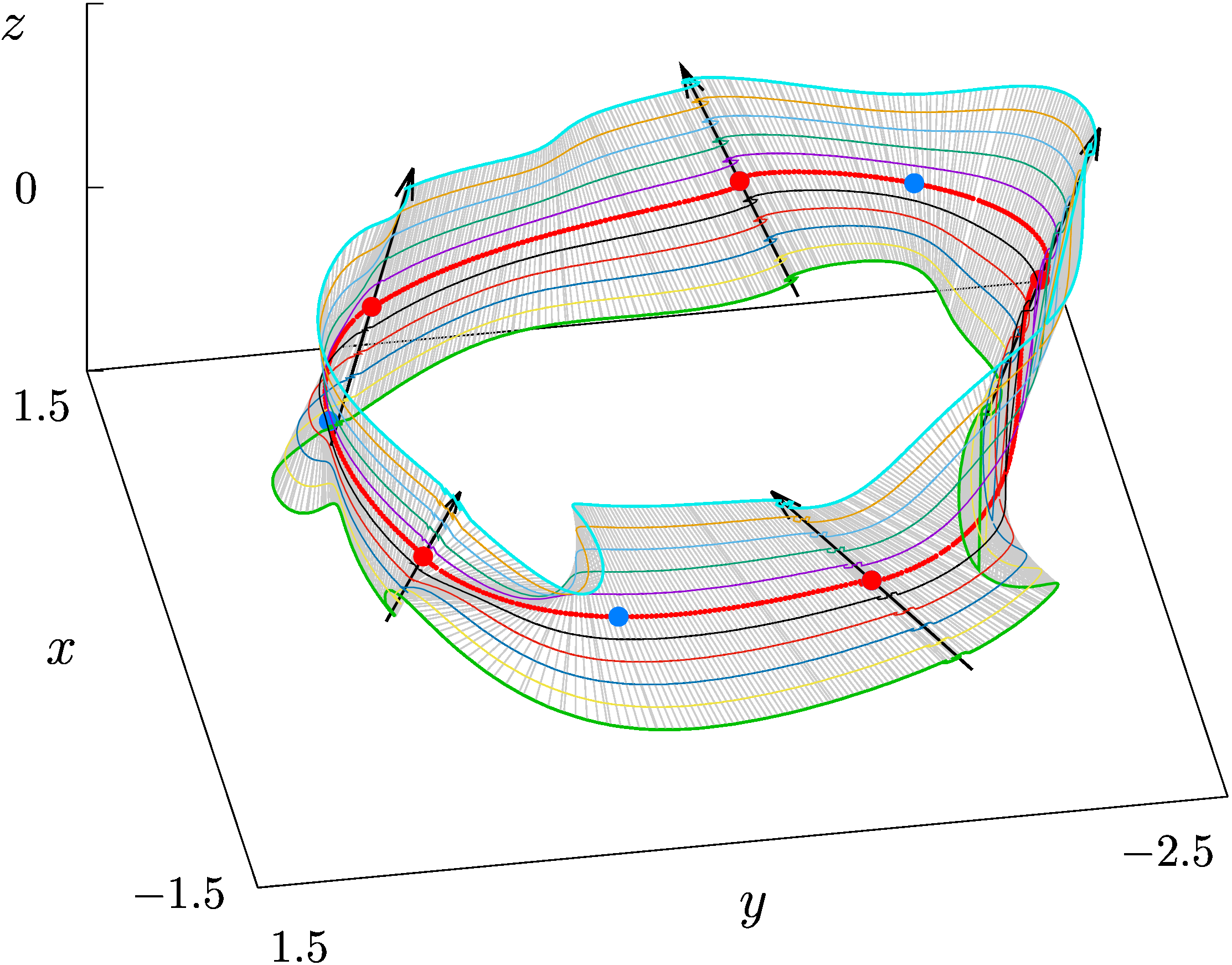}{b}\\
\Includesubgraphics{0.3\textwidth}{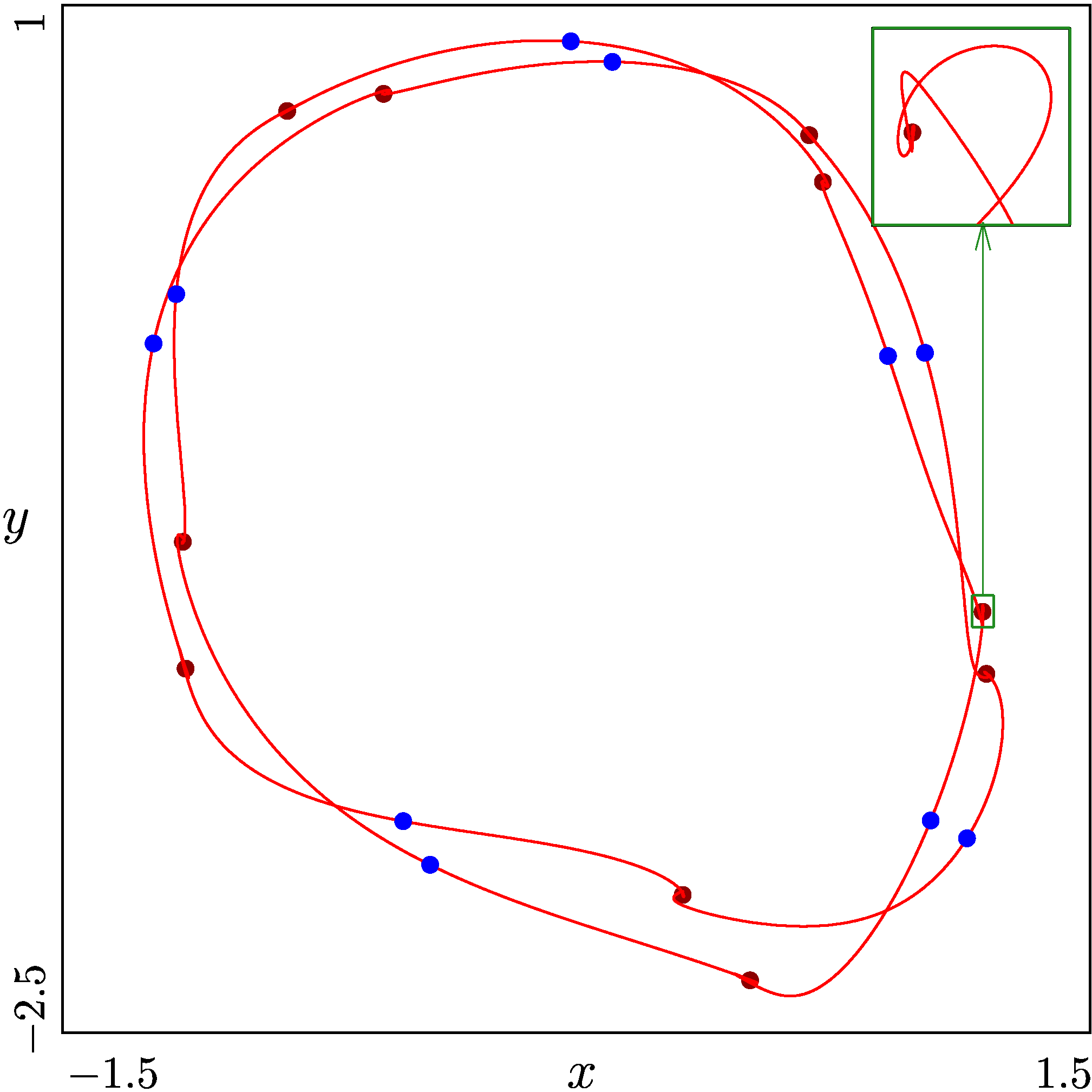}{c}

\caption{\label{fig:loop:SF:doubling}Loop doubling of a saddle-focus ICC in 
map~\eqref{map:cylinder}. (a) The ICC before the bifurcation; (b) Eigenvectors of $J_{f^5}$ corresponding to the doubling 
eigenvalue evaluated at the points of the ICC; (c) The ICC after the bifurcation.
Parameters: $R=1.19$, $E=0.1$,
$\theta=85$, $D=0$,
(a) and (b) $C=-0.999$,
(c) $C=-1.03$.
}
\end{figure*}

Since a torus-doubling bifurcation involves a period doubling of a stable cycle, the real-valued eigenvector (along which the doubling occurs) cannot be a part of the ICC. Therefore, the ICC structure must be formed by a saddle-focus connection. For Fig.~\ref{saddle-focus}(a), since the real eigendirections are along the ICC, it cannot be a period-doubling direction. Therefore, a torus doubling cannot occur for a structure like this.

Fig.~\ref{fig:loop:SF:doubling} shows an ICC with saddle focus connection for the system (\ref{map:cylinder}) that bifurcates to two disjoint loops. The only real eigenvalue of $J_{f^5}$ evaluated at the points
of the stable cycle before the bifurcation is $\approx -0.99501$; the other two eigenvalues are complex. The set of eigenvectors corresponding to the doubling eigenvalue is shown in Fig.~\ref{fig:loop:SF:doubling}(c). Although the direction of the vector changes at different points in the loop, the structure is topologically a cylinder.

\section{Conclusions}

In this work, we have developed a unified and practical framework for
predicting the type of doubling bifurcation (loop doubling or length
doubling) exhibited by invariant closed curves (ICCs) in 3D maps. Our
approach resolves the critical limitations of existing methods by
leveraging the topology of the tangent eigenspace (the eigenvectors
associated with the eigenvalue approaching $-1$ along the
ICC). Specifically, a cylindrical topology predicts loop doubling (two
disjoint ICCs with iterates flipping between them), while a Möbius
strip topology predicts length doubling (a single ICC of double
length). Unlike earlier methods, this technique provides a single,
universally applicable framework for both quasiperiodic and resonant
ICCs. Furthermore, we prove that resonant ICCs with even-period cycles
can only undergo length doubling due to topological constraints,
which explains why loop doubling is observed less frequently than length doubling.

\nonumsection{Acknowledgments} 

\noindent 
S. Banerjee acknowledges ﬁnancial support
from J. C. Bose Grant of the Anusandhan National Research Foundation, Govt. of India, No.
JBR/2020/000049. The work of V. Avrutin and Sushko were supported by the German Research Foundation within the scope of the project ``Global bifurcation phenomena in discontinuous piecewise-smooth maps in theory and applications for power converter systems''. S Biswas acknowledges financial support from the WISE program of DAAD (no. 57698568) to carry out a summer internship in Germany during which some of the work was done.

\bibliographystyle{ws-ijbc}
\bibliography{sb}

\end{document}